\documentclass[11pt,leqno]{article}
\usepackage{amsmath, amssymb, amsfonts, amsthm, amscd}
\input pstricks

\swapnumbers
\theoremstyle{plain}
\newtheorem{theorem}{Theorem}[subsection]
\newtheorem{lemma}[theorem]{Lemma}
\newtheorem{proposition}[theorem]{Proposition}
\newtheorem{corollary}[theorem]{Corollary}

\theoremstyle{definition}
\newtheorem{definition}[theorem]{Definition}
\newtheorem{example}[theorem]{Example}
\newtheorem{examples}[theorem]{Examples}
\newtheorem{remark}[theorem]{Remark}
\newtheorem{remarks}[theorem]{Remarks}

\newcommand{\diag}{\operatorname{diag}}
\newcommand{\id}{\operatorname{id}}
\newcommand{\reg}{{\operatorname{reg}}}
\newcommand{\codim}{\operatorname{codim}}
\newcommand{\supp}{\operatorname{supp}}
\newcommand{\Cl}{\operatorname{Cl}}
\newcommand{\GL}{\operatorname{GL}}
\newcommand{\Gr}{\operatorname{Gr}}
\newcommand{\Pic}{\operatorname{Pic}}

\renewcommand{\div}{\operatorname{div}}

\newcommand{\bC}{{\mathbb C}}
\newcommand{\bP}{{\mathbb P}}
\newcommand{\bQ}{{\mathbb Q}}
\newcommand{\bR}{{\mathbb R}}
\newcommand{\bZ}{{\mathbb Z}}

\newcommand{\cE}{{\mathcal E}}
\newcommand{\cF}{{\mathcal F}}
\newcommand{\cG}{{\mathcal G}}
\newcommand{\cI}{{\mathcal I}}
\newcommand{\cL}{{\mathcal L}}
\newcommand{\cM}{{\mathcal M}}
\newcommand{\cN}{{\mathcal N}}
\newcommand{\cO}{{\mathcal O}}
\newcommand{\cX}{{\mathcal X}}

\newcommand{\uv}{\underline{v}}
\newcommand{\uw}{\underline{w}}

\author{Michel Brion}

\begin{document}

\title{Lectures on the geometry of flag varieties}

\date{}

\maketitle

\section*{Introduction}

In these notes, we present some fundamental results concerning
flag varieties and their Schubert varieties. By a flag variety,
we mean a complex projective algebraic variety $X$, homogeneous under
a complex linear algebraic group. The orbits of a Borel subgroup form
a stratification of $X$ into Schubert cells. These are isomorphic to
affine spaces; their closures in $X$ are the Schubert varieties,
generally singular.

\bigskip

The classes of the Schubert varieties form an additive basis of the
cohomology ring $H^*(X)$, and one easily shows that the structure
constants of $H^*(X)$ in this basis are all non-negative. Our main
goal is to prove a related, but more hidden, statement in the
Grothendieck ring $K(X)$ of coherent sheaves on $X$. This ring admits
an additive basis formed of structure sheaves of Schubert varieties,
and the corresponding structure constants turn out to have alternating
signs.

\bigskip

These structure constants admit combinatorial expressions in the case
of Grassmannians: those of $H^*(X)$ (the Littlewood-Richardson
coefficients) have been known for many years, whereas those of $K(X)$ 
were only recently determined by Buch \cite{Bu02}. This displayed their 
alternation of signs, and Buch conjectured that 
this property extends to all the flag varieties. In this setting, 
the structure constants of the cohomology ring (a fortiori, those 
of the Grothendieck ring) are yet combinatorially elusive, and 
Buch's conjecture was proved in \cite{Br02} by purely algebro-geometric 
methods.

\bigskip

Here we have endeavoured to give a self-contained exposition of 
this proof. The main ingredients are geometric properties of 
Schubert varieties (e.g., their normality), and vanishing theorems for
cohomology of line bundles on these varieties (these are deduced from
the Kawamata-Viehweg theorem, a powerful generalization of the Kodaira
vanishing theorem in complex geometry). Of importance are also
the intersections of Schubert varieties with opposite Schubert
varieties. These ``Richardson varieties'' are systematically used in
these notes to provide geometric explanations for many formulae in the 
cohomology or Grothendieck ring of flag varieties.

\bigskip

The prerequisites are familiarity with algebraic
geometry (for example, the contents of the first three chapters of
Hartshorne's book \cite{Har77}) and with some algebraic topology
(e.g., the book \cite{GrHa81} by Greenberg and Harper). But no
knowledge of algebraic groups is required. In fact, we have presented
all the notations and results in the case of the general linear group,
so that they may be extended readily to arbitrary connected, reductive
algebraic groups by readers familiar with their structure theory.

\bigskip

Thereby, we do not allow ourselves to use the rich algebraic
and combinatorial tools which make Grassmannians and varieties of
complete flags so special among all the flag varieties. For these
developments of Schubert calculus and its generalizations, the reader
may consult the seminal article \cite{LaSc83}, the books \cite{Fu97},
\cite{FuPr98}, \cite{Man98}, and the notes of Buch \cite{Bu04} and 
Tamvakis \cite{Ta04} in this volume. On the other hand, the notes of
Duan in this volume \cite{Du04} provide an introduction to the
differential topology of flag varieties, regarded as homogeneous
spaces under compact Lie groups, with applications to Schubert calculus.

\bigskip

The present text is organized as follows. The first section discusses
Schubert cells and varieties, their classes in the cohomology ring,
and the Picard group of flag varieties. In the second section, we
obtain restrictions on the singularities of Schubert varieties, and
also vanishing theorems for the higher cohomology groups of line
bundles on these varieties. The third section is devoted to a
degeneration of the diagonal of a flag variety into unions of products
of Schubert varieties, with applications to the Grothendieck group. In
the fourth section, we obtain several ``positivity'' results in
this group, including a solution of Buch's conjecture. Each section
begins with a brief overview of its contents, and ends with
bibliographical notes and open problems.

\bigskip

These notes grew out of courses at the Institut Fourier (Grenoble) in
the spring of 2003, and at the mini-school ``Schubert Varieties''
(Banach Center, Warsaw) in May 2003. I am grateful to the organizers of
this school, Piotr Pragacz and Andrzej Weber, for their invitation and
encouragements. I also thank the auditors of both courses, especially
Dima Timashev, for their attention and comments. 

\bigskip

\noindent
{\bf Conventions.} 

Throughout these notes, we consider algebraic varieties over the field
$\bC$ of complex numbers. We follow the notation and terminology of
\cite{Har77}; in particular, varieties are assumed to be irreducible. 
Unless otherwise stated, subvarieties are assumed to be closed.

\vfill\eject

\section{Grassmannians and flag varieties}

We begin this section by reviewing the definitions and fundamental
properties of Schubert varieties in Grassmannians and varieties
of complete flags. Then we introduce the Schubert classes in the
cohomology ring of flag varieties, and we study their multiplicative 
properties. Finally, we describe the Picard group of flag varieties, 
first in terms of Schubert divisors, and then in terms of homogeneous
line bundles; we also sketch the relation of the latter to 
representation theory.
 
\subsection{Grassmannians}

The {\it Grassmannian} $\Gr(d,n)$ is the set of 
$d$-dimensional linear subspaces of $\bC^n$. Given such a subspace $E$
and a basis $(v_1,\ldots,v_d)$ of $E$, the exterior product 
$v_1\wedge \cdots \wedge v_d \in \bigwedge^d \bC^n$ 
only depends on $E$ up to a non-zero scalar multiple. In other words,
the point  
$$
\iota(E) := \left[ v_1\wedge \cdots \wedge v_d \right]
$$
of the projective space $\bP(\bigwedge^d \bC^n)$ only depends on $E$.
Further, $\iota(E)$ uniquely determines $E$, so that the map $\iota$ 
identifies $\Gr(d,n)$ with the image in $\bP(\bigwedge^d \bC^n)$ of
the cone of decomposable $d$-vectors in $\bigwedge^d \bC^n$. It
follows that $\Gr(d,n)$ is a subvariety of the projective space
$\bP(\bigwedge^d \bC^n)$; the map
$$
\iota:\Gr(d,n) \to \bP(\bigwedge^d \bC^n)
$$ 
is the \emph{Pl\"ucker embedding}.

The general linear group 
$$
G:=GL_n(\bC)
$$ 
acts on the variety
$$
X:=\Gr(d,n)
$$ 
via its natural action on $\bC^n$. Clearly, $X$ is a unique $G$-orbit,
and the Pl\"ucker embedding is equivariant with respect to the action
of $G$ on $\bP(\bigwedge^d \bC^n)$ arising from its linear action
on $\bigwedge^d \bC^n$. Let $\left(e_1,\ldots,e_n\right)$
denote the standard basis of $\bC^n$, then the isotropy group of the
subspace $\langle e_1,\ldots,e_d\rangle$ is
$$
P :=\left\{
\begin{pmatrix}
a_{1,1} & \ldots & a_{1,d} & a_{1,d+1} & \ldots & a_{1,n}\cr
\vdots & \ddots & \vdots & \vdots & \ddots & \vdots\cr
a_{d,1} & \ldots & a_{d,d} & a_{d,d+1} & \ldots & a_{d,n}\cr
0 & \ldots & 0 & a_{d+1,d+1} & \ldots & a_{d+1,n}\cr
\vdots & \ddots & \vdots & \vdots & \ddots & \vdots\cr
0 & \ldots & 0 & a_{n,d+1} & \ldots & a_{n,n}\cr
\end{pmatrix}
\right\}
$$
(this is a {\it maximal parabolic subgroup} of $G$). Thus, $X$ is the
homogeneous space $G/P$. As a consequence, the algebraic variety $X$
is nonsingular, of dimension $\dim(G) - \dim(P) = d(n-d)$.

For any multi-index $I : = (i_1,\ldots,i_d )$, where 
$1\leq i_1 < \ldots < i_d \leq n$, we denote by $E_I$ the
corresponding {\it coordinate subspace} of $\bC ^n$, i.e., 
$E_I=\langle e_{i_1},\ldots, e_{i_d} \rangle \in X$. In particular, 
$E_{1,2,\ldots,d}$ is the {\it standard coordinate subspace}
$\langle e_1,\ldots,e_d\rangle$. We may now state the following
result, whose proof is straightforward.

\begin{proposition}\label{grass1}
(i) The $E_I$ are precisely the $T$-fixed points in $X$, where
$$T := \left\{
\begin{pmatrix}
a_{1,1} & 0 & \ldots & 0\cr
0 & a_{2,2} & \ldots & 0\cr
\vdots & \vdots & \ddots & \vdots\cr
0 & 0 & \ldots & a_{n,n}\cr
\end{pmatrix}
\right\}\subseteq GL_n(\bC)
$$
is the subgroup of diagonal matrices (this is a {\rm maximal torus} of
$G$). 

\noindent
(ii) $X$ is the disjoint union of the orbits $B E_I$, where
$$
B := \left\{
\begin{pmatrix}
a_{1,1} & a_{1,2} & \ldots & a_{1,n}\cr
0 & a_{2,2} & \ldots & a_{2,n}\cr
\vdots & \vdots & \ddots & \vdots\cr
0 & 0 & \ldots & a_{n,n}\cr
\end{pmatrix}
\right\}\subseteq GL_n(\bC)
$$
is the subgroup of upper triangular matrices (this is a {\rm Borel
subgroup} of $G$). 
\end{proposition}

\begin{definition}
The {\it Schubert cells} in the Grassmannian are the orbits
$C_I := B E_I$, i.e., the $B$-orbits in $X$. The closure in
$X$ of the Schubert cell $C_I$ (for the Zariski topology) is called
the {\it Schubert variety} $X_I := \overline {C_I}$.
\end{definition}

Note that $B$ is the semi-direct product of $T$ with the normal
subgroup
$$
U := \left\{
\begin{pmatrix}
1 & a_{1,2} & \ldots & a_{1,n}\cr
0 & 1 & \ldots & a_{2,n}\cr
\vdots & \vdots & \ddots & \vdots\cr
0 & 0 & \ldots & 1\cr
\end{pmatrix}
\right\}
$$
(this is a {\it maximal unipotent subgroup} of $G$). Thus, we also
have $C_I = U E_I$ : the Schubert cells are just the $U$-orbits in $X$.

Also, the isotropy group $U_{E_I}$ is the subgroup of $U$ where
$a_{ij}=0$ whenever $i\notin I$ and $j\in I$. Let $U^I$ be the
``complementary'' subset of $U$, defined by $a_{ij}=0$ if $i\in I$
or $j\notin I$. Then one checks that $U^I$ is a subgroup of $U$, and
the map $U^I\to X$, $g\mapsto gE_I$ is a locally closed embedding with
image $C_I$. It follows that $C_I$ is a locally closed subvariety of
$X$, isomorphic to the affine space $\bC^{\vert I \vert}$, where 
$\vert I \vert :=\sum_{j=1}^d (i_j - j)$.
Thus, its closure $X_I$ is a projective variety of dimension 
$\vert I \vert$.

Next we present a geometric characterization of Schubert cells and
varieties (see e.g. \cite{Fu97} 9.4).

\begin{proposition}\label{grass2}
(i) $C_I$ is the set of $d$-dimensional subspaces $E\subset \bC^n$
such that 
$$
\dim(E\cap \langle e_1,\ldots,e_j \rangle )
=\# \left\{ k ~\vert~ 1\le k\le d, ~i_k<j\right\}, \quad 
\text{for}\quad j=1,\ldots,n.
$$

\noindent
(ii) $X_I$ is the set of $d$-dimensional subspaces $E\subset \bC^n$
such that  
$$
\dim(E\cap \langle e_1,\ldots,e_j \rangle ) \geq 
\# \left\{ k ~\vert~ 1\le k\le d,~ i_k<j\right\}, \quad 
\text{for} \quad j=1,\ldots,n.
$$
Thus, we have 
$$
X_I = \bigcup_{J\leq I} C_J,
$$ 
where $J\leq I$ if and only if $j_k\leq i_k$ for all $k$.
\end{proposition}

\begin{examples}\label{gras} 
1) For $d=1$, the Grassmannian is just the projective
space $\bP^{n-1}$, and the Schubert varieties form a flag of linear
subspaces $X_0\subset X_1\subset\cdots\subset X_n$, where 
$X_j \cong \bP^{j-1}$.  

\medskip

\noindent
2) For $d=2$ and $n=4$ one gets the following poset of Schubert varieties:
$$
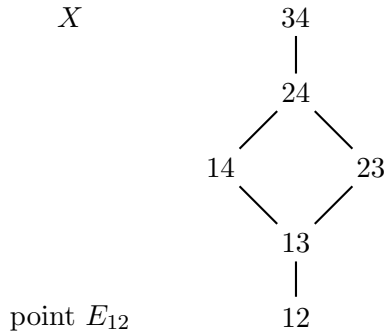
\pspicture(-5,-3)(5,3)
\rput(-3,2){$X$}
\rput(-3,-2){point $E_{12}$}
\rput(0,2){34}
\psline(0,1.75)(0,1.25)
\rput(0,1){24}
\psline(-.25,.75)(-.75,.25)
\psline(.25,.75)(.75,.25)
\rput(1,0){23}
\rput(-1,0){14}
\rput(0,-2){12}
\psline(0,-1.72)(0,-1.25)
\rput(0,-1){13}
\psline(.25,-.75)(.75,-.25)
\psline(-.25,-.75)(-.75,-.25)
\endpspicture
$$
Further, the Schubert variety $X_{24}$ is singular. Indeed, one checks
that $X\subset \bP(\bigwedge^2\bC^4)=\bP^5$ is defined by one quadratic
equation (the Pl\"ucker relation). Further, $X_{24}$ is the
intersection of $X$ with its tangent space at the point $E_{12}$. 
Thus, $X_{24}$ is a quadratic cone with vertex $E_{12}$, its
unique singular point.

\medskip

\noindent
3) For arbitrary $d$ and $n$, the Schubert variety
$X_{1,2,\ldots,d}$ is just the point $E_{1,2,\ldots,d}$, whereas 
$X_{n-d+1,n-d+2,\ldots,n}$ is the whole Grassman\-nian. On the other hand, 
$X_{n-d,n-d+2,\dots,n}$ consists of those $d$-dimensional subspaces
$E$ that meet $\langle e_1,\dots,e_{n-d}\rangle$: it is the intersection
of $X$ with the hyperplane of $\bP(\bigwedge^d \bC^n)$ where the 
coordinate on $e_{n-d+1}\wedge\cdots\wedge e_n$ vanishes. 

Since $X$ is the disjoint union of the open Schubert cell 
$C_{n-d+1,n-d+2,\ldots,n}\cong\bC^{d(n-d)}$ with the irreducible
divisor $D := X_{n-d,n-d+2,\ldots,n}$, any divisor in $X$ is linearly
equivalent to a unique integer multiple of $D$. Equivalently, any line
bundle on $X$ is isomorphic to a unique tensor power of the line
bundle $L := \cO_X(D)$, the pull-back of $\cO(1)$ via the Pl\"ucker
embedding. Thus, {\it the Picard group $\Pic(X)$ is freely generated
by the class of the very ample line bundle $L$.}
\end{examples}

We may re-index Schubert varieties in two ways:

\noindent
{\it 1. By partitions:} with any multi-index $I=(i_1,\ldots,i_d)$
we associate the partition $\lambda=(\lambda_1,\ldots,\lambda_d)$,
where $\lambda_j := i_j - j$ for $j=1,\ldots,d$. We then write
$X_{\lambda}$ instead of $X_I$.

This yields a bijection between the set of multi-indices
$I=(i_1,\ldots,i_d)$ such that $1\leq i_1<\ldots<i_d\leq n$, and the
set of tuples of integers $\lambda=(\lambda_1,\dots,\lambda_d)$
satisfying $0\leq\lambda_1\leq\ldots\leq\lambda_d\leq n-d$. This is
the set of {\it partitions with $\le d$ parts of size $\le n-d$}.

The {\it area} of the partition $\lambda$ is the number 
$|\lambda |:=\sum_{j=1}^{d}\lambda_j=|I|$.
With this indexing, the dimension of $X_\lambda$ is the area of
$\lambda$; further, $X_{\mu}\subseteq X_{\lambda}$ if and only if
$\mu \le \lambda$, that is, $\mu_j\le \lambda_j$ for all $j$.

Alterrnatively, one may associate with any multi-index
$I=(i_1,\ldots,i_d)$ the {\it dual partition} 
$(n-i_d,n-1-i_{d-1},\ldots,n-d+1-i_1)$. This is still a partition with
$\le d$ parts of size $\le n-d$, but now its area is the codimension of
the corresponding Schubert variety. This indexing is used in the
notes of Buch \cite{Bu04} and Tamvakis \cite{Ta04}.

\noindent
{\it 2. By permutations:} with a multi-index $I=(i_1,\ldots,i_d)$
we associate the permutation $w$ of the set $\{1,2,\ldots,n\}$,
defined as follows: $w(k)=i_k$ for $k=1,\ldots,d$, whereas 
$w(d+k)$ is the $k$-th element of the ordered set
$\{1,\ldots,n\}\setminus I$ for $k=1,\ldots,n-d$.
This sets up a bijection between the multi-indices and the
permutations $w$ such that 
$w(1)< w(2) < \cdots < w(d)$ and $w(d+1) < \cdots < w(n)$.
These permutations form a system of representatives of the
coset space $S_n/(S_d\times S_{n-d})$, where $S_n$ denotes the
permutation group of the set $\{1,2,\ldots,n\}$, and 
$S_d\times S_{n-d}$ is its subgroup stabilizing the subset
$\{1,2,\ldots,d\}$ (and $\{d+1,d+2,\ldots,n\}$). Thus, we may
parametrize the $T$-fixed points of $X$, and hence the Schubert
varieties, by the map 
$S_n/(S_d\times S_{n-d}) \rightarrow X$,
$w(S_d\times S_{n-d})\mapsto E_{w(1),\ldots,w(d)}$.
This parametrization will be generalized to all flag varieties 
in the next subsection.

\subsection{Flag varieties}

\medskip

Given a sequence $(d_1,\ldots,d_m)$ of positive integers with sum $n$, a
{\it flag of type $(d_1,\ldots,d_m)$ in $\bC^n$} 
is an increasing sequence of linear subspaces 
$$
0=V_0\subset V_1\subset V_2\subset\ldots\subset V_m=\bC^n
$$ 
such that $\dim({V_j/V_{j-1}})=d_j$ for $j=1,\ldots,m$.
The {\it coordinate flags} are those consisting of coordinate subspaces.

Let $X(d_1,\ldots,d_m)$ denote the set of flags of type
$(d_1,\ldots,d_m)$. For example, $X(d,n-d)$ is just the Grassmannian
$\Gr(d,n)$. More generally, $X(d_1,\ldots,d_m)$ is a subvariety of
the product of the Grassmannians $\Gr(d_i,n)$, called the 
{\it partial flag variety} of type $(d_1,\ldots, d_m)$.

The group $G=\GL_n(\bC)$ acts transitively on $X(d_1,\ldots,d_m)$.
Let $P=P(d_1,\ldots,d_m)$ be the isotropy group of the 
{\it standard flag} (consisting of the standard coordinate
subspaces). Then $P(d_1,\ldots,d_m)$ consists of the block upper 
triangular invertible matrices with diagonal blocks of sizes
$d_1,\ldots,d_m$. In particular, $P(d_1,\ldots,d_m)$ contains $B$; in
fact, all subgroups of $G$ containing $B$ occur in this way.
(These subgroups are the {\it standard parabolic subgroups of $G$}).
Since $X \cong G/P$, it follows that $X$ is nonsingular of dimension
$\sum_{1\le i<j\le m}d_i d_j$.

In particular, we have the variety $X:=X(1,\ldots,1)$ of 
{\it complete flags}, also called the {\it full flag variety}; 
it is the homogeneous space $G/B$, of dimension $n(n-1)/2$. 
By sending any complete flag to the corresponding partial flag of 
a given type $(d_1,\ldots,d_m)$, we obtain a morphism
$$
f:X=G/B\to G/P(d_1,\ldots,d_m)=X(d_1,\ldots,d_m).
$$
Clearly, $f$ is $G$-equivariant with fiber $P/B$ at the base point
$B/B$ (the standard complete flag). 
Thus, $f$ is a fibration with fiber being the product of varieties of
complete flags in $\bC^{d_1}$, $\ldots$, $\bC^{d_m}$. This allows us
to reduce many questions regarding flag varieties to the case of the
variety of complete flags; see Example \ref{flag} below for details on
this reduction. Therefore, we will mostly concentrate on the full flag
variety. 

We now introduce Schubert cells and varieties in $G/B$.
Observe that the complete coordinate flags correspond to the
permutations of the set $\{1,\ldots,n\}$, by assigning to the flag
$$
0 \subset \langle e_{i_1}\rangle \subset \cdots
\subset \langle e_{i_1},e_{i_2},\ldots, e_{i_k}\rangle 
\subset\cdots
$$
the permutation $w$ such that $w(k)= i_k$ for all $k$.  We
regard the permutation group $S_n$ as a subgroup of $\GL_n(\bC)$ via 
its natural action on the standard basis $(e_1,\ldots,e_n)$. Then the
(complete) coordinate flags are exactly the $F_w := w F$, where 
$F$ denotes the standard complete flag. Further, $S_n$ identifies to
the quotient $W:=N_G(T)/T$, where $N_G(T)$ denotes the normalizer of
$T$ in $G$. (In other words, $S_n$ is the {\it Weyl group of $G$ with
respect to $T$}). 

We may now formulate an analogue of Proposition \ref{grass1} (see
e.g. \cite{Fu97} 10.2 for a proof).

\begin{proposition} (i) The fixed points of $T$ in $X$ are the
coordinate flags $F_w$, $w\in W$.

\noindent
(ii) $X$ is the disjoint union of the orbits 
$C_w:=B F_w=U F_w$, where $w\in W$.

\noindent
(iii) Let $X_w:=\overline {C_w}$ (closure in the Zariski
topology of $X$), then 
$$
X_w=\bigcup_{v\in W, \; v \leq w} C_v,
$$ 
where $v \leq w $ if and only if we have 
$(v(1),\dots,v(d))_{\hbox{\rm r.t.i.v.}}
\leq(w(1),\dots,w(d))_{\hbox{\rm r.t.i.v.}}$ 
for $d=1,\ldots,n-1$ (here r.t.i.v. stands for ``reordered to
increasing values'').
\end{proposition}

\begin{definition}
$C_w := B F_w$ is a {\it Schubert cell}, and
$X_w := \overline{C_w}$ is the corresponding 
{\it Schubert variety}. The partial ordering $\le$ on $W$ is the 
{\it Bruhat order}.
\end{definition}

By the preceding proposition, we have $X_v\subseteq X_w$
if and only if this holds for the images of $X_v$ and $X_w$ in
$\Gr(d,n)$, where $d=1,\ldots,n-1$. Together with Proposition
\ref{grass2}, this yields a geometric characterization of   
the Bruhat order on Schubert varieties. Also, note that the $T$-fixed
points in $X_w$ are the coordinate flags $F_v$, where 
$v\in W$ and $v \le w$.

We now describe the Schubert cells $U F_w$. Note that the
isotropy group 
$$
U_{F_w} = U\cap w Uw^{-1}=:U_w
$$ 
is defined by $a_{i,j}=0$ whenever $i<j$ and $w^{-1}(i)<w^{-1}(j)$. 
Let $U^w$ be the ``complementary'' subset of $U$, defined by 
$a_{ij}=0$ whenever $i<j$ and $w^{-1}(i)>w^{-1}(j)$. Then 
$U^w = U\cap w U^-w^{-1}$ is a subgroup, and one checks that the
product map $U^w\times U_w\to U$ is an isomorphism of varieties. 
Hence the map $U^w \to C_w$, $g\mapsto g F_w$ is an isomorphism as
well. It follows that each $C_w$ is an affine space of dimension 
$$
\# \{ (i,j)~\vert~1\le i<j\le n, ~w^{-1}(i)>w^{-1}(j)\}
= \# \{ (i,j)~\vert~1\le i<j\le n, ~w(i)>w(j)\}.
$$
The latter set consists of the {\it inversions} of the permutation
$w$; its cardinality is the {\it length} of $w$, denoted by
$\ell(w)$. Thus, $C_w\cong \bC^{\ell(w)}$.

More generally, we may define Schubert cells and varieties in any
partial flag variety $X(d_1,\ldots,d_m) =G/P$, where
$P=P(d_1,\ldots,P_m)$; these are parametrized by the coset space 
$S_n/(S_{d_1}\times\cdots\times S_{d_m}) =: W/W_P$. 

Specifically, each right coset mod $W_P$ contains a unique permutation
$w$ such that we have $w(1) < \cdots < w(d_1)$, 
$w(d_1+1) < \cdots < w(d_1+d_2)$, $\ldots$, 
$w(d_1+\cdots + d_{m-1}+1)< \cdots < w(d_1 + \cdots + d_m) = w(n)$.
Equivalently, $w \le wv$ for all $v\in W_P$. This defines the set
$W^P$ of {\it minimal representatives} of $W/W_P$. 

Now the Schubert cells in $G/P$ are the orbits 
$C_{wP}:=BwP/P = UwP/P \subset G/P$ ($w\in W^P$), and the Schubert
varieties $X_{wP}$ are their closures. One checks that the map
$f:G/B\to G/P$ restricts to an isomorphism 
$C_w = BwB/B \cong BwP/P= C_{wP}$, and hence to a birational morphism
$X_w \to X_{wP}$ for any $w\in W^P$.

\begin{examples}\label{flag} 
1) The Bruhat order on $S_2$ is just
$$\pspicture(-1,-1)(1,1)
\rput(0,.5){$(21)$}
\psline(0,.25)(0,-.25)
\rput(0,-.5){$(12)$}
\endpspicture$$

The picture of the Bruhat order on $S_3$ is 
$$
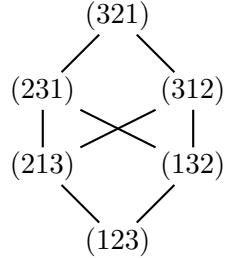
\pspicture(-2,-2)(2,2)
\rput(0,1.5){$(321)$}
\psline(-.25,1.25)(-.75,.75)
\psline(.25,1.25)(.75,.75)
\rput(1,.5){$(312)$}
\rput(-1,.5){$(231)$}
\psline(-1,.25)(-1,-.25)
\psline(.5,.25)(-.5,-.25)
\psline(-.5,.25)(.5,-.25)
\psline(1,.25)(1,-.25)
\rput(1,-.5){$(132)$}
\rput(-1,-.5){$(213)$}
\psline(-.25,-1.25)(-.75,-.75)
\psline(.25,-1.25)(.75,-.75)
\rput(0,-1.5){$(123)$}
\endpspicture$$

\medskip

\noindent
2) Let $w_o := (n,n-1,\dots,1)$, the order-reversing permutation. Then
$X= X_{w_o}$, i.e., $w_o$ is the unique maximal element of the Bruhat
order on $W$. Note that $w_o^2=\id$, and $\ell(w_o w) = \ell(w_o)-\ell(w)$ 
for any $w\in W$. 

\medskip

\noindent
3) The permutations of length $1$ are exactly the 
{\it elementary transpositions} $s_1,\ldots,s_{n-1}$, where each $s_i$
exchanges the indices $i$ and $i+1$ and fixes all other indices.
The corresponding Schubert varieties are the {\it Schubert curves}
$X_{s_1},\ldots,X_{s_{n-1}}$. In fact, $X_{s_i}$ may be identified
with the set of $i$-dimensional subspaces $E\subset \bC^n$ such that  
$$
\langle e_1,\ldots,e_{i-1}\rangle \subset E \subset 
\langle e_1,\ldots,e_{i+1}\rangle.
$$ 
Thus, $X_{s_i}$ is the projectivization of the quotient space 
$\langle e_1,\ldots,e_{i+1}\rangle/\langle e_1,\ldots,e_{i-1}\rangle$,
so that $X_{s_i}\cong\bP^1$.

\medskip

\noindent
4) Likewise, the Schubert varieties of codimension $1$ are 
$X_{w_o s_1},\ldots,X_{w_o s_{n-1}}$, also called the 
{\it Schubert divisors}.

\medskip

\noindent
5) Apart from the Grassmannians, the simplest partial flag variety is
the {\it incidence variety} $I=I_n$ consisting of the pairs
$(V_1,V_{n-1})$, where $V_1\subset \bC^n$ is a line, and
$V_{n-1}\subset \bC^n$ is a hyperplane containing $V_1$. Denote by
$\bP^{n-1}=\bP(\bC^n)$ (resp.~$\check\bP^{n-1}=\bP((\bC^n)^*)$) the
projective space of lines (resp.~hyperplanes) in $\bC^n$, then 
$I\subset \bP^{n-1}\times \check\bP^{n-1}$
is defined by the bi-homogeneous equation
$$
x_1 y_1 + \cdots + x_n y_n=0,
$$
where $x_1,\ldots,x_n$ are the standard coordinates on $\bC^n$, and 
$y_1,\ldots,y_n$ are the dual coordinates on $(\bC^n)^*$.

One checks that the Schubert varieties in $I$ are the 
$$
I_{i,j}:= \{(V_1,V_{n-1})\in I  \quad \vert \quad 
V_1\subseteq E_{1,\ldots,i} \quad \text{and} \quad 
E_{1,\ldots,j-1}\subseteq V_{n-1}\},
$$
where $1\le i,j\le n$ and $i\ne j$. Thus, $I_{i,j}\subseteq I$ is
defined by the equations 
$$
x_{i+1}= \cdots = x_n = y_1 = \cdots = y_{j-1} = 0.
$$
It follows that {\it $I_{i,j}$ is singular for $1<j<i<n$ with singular 
locus $I_{j-1,i+1}$, and is nonsingular otherwise.} 

\medskip

\noindent
6) For any partial flag variety $G/P$ and any $w\in W^P$, the
pull-back of the Schubert variety $X_{wP}$ under $f:G/B\to G/P$ is
easily seen to be the Schubert variety $X_{ww_{0,P}}$, where $w_{0,P}$
denotes the maximal element of $W_P$. Specifically, if
$P=P(d_1,\ldots,d_m)$ so that 
$W_P = S_{d_1} \times \cdots \times S_{d_m}$, then 
$w_{0,P}=(w_{0,d_1},\ldots,w_{0,d_m})$ with obvious notation.
The products $ww_{0,P}$, where $w\in W^P$, are the 
{\it maximal representatives} of the cosets modulo $W_P$. Thus, $f$
restricts to a locally trivial fibration $X_{ww_{0,P}}\to X_{wP}$
with fiber $P/B$. 

In particular, the preceding example yields many singular Schubert
varieties in the variety of complete flags, by pull-back from the
incidence variety.
\end{examples}

\begin{definition}
The {\it opposite Schubert cell (resp.~variety}) associated with
$w \in W$ is $C^w:= w_o C_{w_o w}$ (resp.
$X^w:= w_o X_{w_o w}$).
\end{definition}

Observe that $C^w=B^- F_w$, where 
$$
B^- := \left\{
\begin{pmatrix}
a_{1,1} & 0 & \ldots & 0 \cr
a_{2,1} & a_{2,2} & \ldots & 0  \cr
\vdots & \vdots & \ddots & \vdots \cr
a_{n,1} & a_{n,2} & \ldots & a_{n,n}  \cr
\end{pmatrix}\right\}=w_o B w_o
$$
(this is the {\it opposite Borel subgroup} to $B$ containing the
maximal torus $T$). Also, $X^w$ has codimension $\ell(w)$ in $X$.

For example, $C^{\id}\cong U^-$ via the map $U^-\to X$, 
$g\mapsto g F$, where $U^-:=w_o U w_o$.
Further, this map is an open immersion. Since $X=G/B$, this is
equivalent to the fact that the product map $U^-\times B\to G$ is an
open immersion (which, of course, may be checked directly).
It follows that the quotient $q: G\to G/B$, $g\mapsto gB$, is a
trivial fibration over $C^{\id}$; thus, by $G$-equivariance, $q$ is 
locally trivial for the Zariski topology. This also holds for any
partial flag variety $G/P$ with the same proof. Likewise, the map
$f:G/B \to G/P$ is a locally trivial fibration with fiber $P/B$.

\subsection{Schubert classes}

This subsection is devoted to the cohomology ring of the full flag
variety. We begin by recalling some basic facts on the homology and
cohomology of algebraic varieties, referring for details to
\cite{Fu97} Appendix B or \cite{FuPr98} Appendix A. We will
consider (co)homology groups with integer coefficients.

Let $X$ be a projective nonsingular algebraic variety of dimension
$n$. Then $X$ (viewed as a compact differentiable manifold of
dimension $2n$) admits a canonical orientation, hence a canonical
generator of the homology group $H_{2n}(X)$ : the 
{\it fundamental class} $[X]$. By Poincar\'e duality, the map  
$H^j(X)\to H_{2n-j}(X)$, $\alpha \mapsto \alpha \cap [X]$ is an
isomorphism for all $j$. 

Likewise, any nonsingular subvariety $Y\subseteq X$ of dimension $p$
has a fundamental class in $H_{2p}(Y)$. Using Poincar\'e duality, the
image of this class in $H_{2p}(X)$ yields the {\it fundamental class}
$[Y]\in H^{2c}(X)$, where $c =n - p$ is the codimension of $Y$. 
In particular, we obtain the fundamental class of a point $[x]$, which
is independent of $x$ and generates the group $H^{2n}(X)$. More
generally, one defines the fundamental class $[Y] \in H^{2c}(X)$ for
any (possibly singular) subvariety $Y$ of codimension $c$.

Given $\alpha$, $\beta$ in the cohomology ring $H^*(X)$, let 
$\langle\alpha,\beta \rangle$ denote the coefficient of thc class
$[x]$ in the cup product $\alpha\cup\beta$. Then $\langle,\rangle$ is
a bilinear form on $H^*(X)$ called the {\it Poincar\'e duality pairing}.
It is non-degenerate over the rationals, even over the integers
in the case where the group $H^*(X)$ is torsion-free.

For any two subvarieties $Y$, $Z$ of $X$, each irreducible component
$C$ of $Y\cap Z$ satisfies $\dim(C)\ge \dim(Y)+\dim(Z)$, i.e.,
$\codim(C)\le \codim(Y) + \codim(Z)$. We say that 
{\it $Y$ and $Z$ meet properly in $X$}, if 
$\codim(C) = \codim(Y) + \codim(Z)$ for each $C$. Then we have in 
$H^*(X)$:
$$
[Y]\cup [Z] = \sum_C m_C \; [C],
$$
where the sum is over all irreducible components of $Y\cap Z$, and
$m_C$ is the {\it intersection multiplicity of $Y$ and $Z$ along $C$},
a positive integer. Further, $m_C=1$ if and only if 
{\it $Y$ and $Z$ meet transversally along $C$}, i.e., there exists
a point $x\in C$ such that: $x$ is a nonsingular point of $Y$ and $Z$,
and the tangent spaces at $x$ satisfy $T_x Y + T_x Z = T_x X$.
Then $x$ is a nonsingular point of $C$, and $T_x C = T_x Y \cap T_x Z$. 

In particular, if $Y$ and $Z$ are subvarieties such that
$\dim(Y) + \dim(Z)=\dim(X)$, then $Y$ meets $Z$ properly if and only
if their intersection is finite. In this case, we have 
$\langle [Y],[Z] \rangle = \sum_{x\in Y\cap Z} m_x$,
where $m_x$ denotes the intersection multiplicity of $Y$ and $Z$ at
$x$. In the case of transversal intersection, this simplifies to 
$\langle [Y],[Z] \rangle = \#(Y\cap Z)$.

Returning to the case where $X$ is a flag variety, we have the
cohomology classes of the Schubert subvarieties, called the 
{\it Schubert classes}. Since $X$ is the disjoint union of the
Schubert cells, {\it the Schubert classes form an additive basis of
$H^*(X)$}; in particular, this group is torsion-free. 

To study the cup product of Schubert classes, we will need a version
of Kleiman's transversality theorem, see \cite{Kl74} or \cite{Har77}
Theorem III.10.8.

\begin{lemma}\label{klei}
Let $Y$, $Z$ be subvarieties of a flag variety $X$
and let $Y_0\subseteq Y$ (resp. $Z_0\subseteq Z$) be nonempty open
subsets consisting of nonsingular points. Then there
exists a nonempty open subset $\Omega$ of $G$ such that: for any 
$g\in \Omega$, $Y$ meets $gZ$ properly, and $Y_0\cap gZ_0$ is
nonsingular and dense in $Y\cap gZ$. Thus, 
$[Y]\cup [Z] = [Y\cap gZ]$ for all $g\in\Omega$.

In particular, if $\dim(Y) + \dim(Z) = \dim(X)$, then $Y$ meets $gZ$
transversally {\rm for general} $g\in G$, that is, for all $g$ in a
nonempty open subset $\Omega$ of $G$. Thus, $Y\cap gZ$ is finite and
$\langle [Y],[Z] \rangle= \#(Y\cap gZ)$, for general $g\in G$. 
\end{lemma}

\begin{proof}
Consider the map $m:G\times Z\to X$, $(g,z)\mapsto gz$. This is a
surjective morphism, equivariant for the action of $G$ on $G\times Z$
by left multiplication on the first factor. Since $X=G/P$, it follows
that $m$ is a locally trivial fibration for the Zariski
topology. Thus, its scheme-theoretic fibers are varieties of dimension
$\dim(G)+\dim(Z)-\dim(X)$.

Next consider the fibered product $V:=(G\times Z)\times_X Y$ and the
pull-back $\mu:V\to Y$ of $m$. Then $\mu$ is also a locally trivial
fibration with fibers being varieties. It follows that the scheme $V$
is a variety of dimension $\dim(G)+\dim(Z)-\dim(X)+\dim(Y)$.

Let $\pi:V\to G$ be the composition of the projections
$(G\times Z)\times_X Y\to G\times Z \to G$. Then the fiber of $\pi$ at
any $g\in G$ may be identified with the scheme-theoretic intersection
$Y\cap gZ$. Further, there exists a nonempty open subset $\Omega$ of
$G$ such that the fibers of $\pi$ at points of $\Omega$ are either
empty or equidimensional of dimension $\dim(Y)+\dim(Z)-\dim(X)$, i.e.,
of codimension $\codim(Y)+\codim(Z)$. This shows that $Y$ meets $gZ$
properly for any $g\in\Omega$. 

Likewise, the restriction $m_0:G\times Z_0\to X$ is a locally
trivial fibration with nonsingular fibers, so that the fibered product 
$V_0:=(G\times Z_0)\times_X Y_0$ is a nonempty open subset of $V$
consisting of nonsingular points. By generic smoothness, it follows
that $Y_0\cap gZ_0$ is nonsingular and dense in $Y\cap gZ$, for all
$g$ in a (possibly smaller) nonempty open subset of $G$. This implies,
in turn, that all intersection multiplicities of $Y\cap Z$ are $1$.

Thus, we have $[Y] \cup [gZ] = [Y\cap gZ]$ for any $g\in \Omega$.
Further, $[Z]=[gZ]$ as $G$ is connected, so that 
$[Y] \cup [Z] = [Y\cap gZ]$.
\end{proof}

As a consequence, in the full flag variety $X$, any Schubert variety
$X_w$ meets properly any opposite Schubert variety $X^v$. (Indeed, the
open subset $\Omega$ meets the open subset 
$BB^- = BU^- \cong B\times U^-$ of $G$; further, 
$X_w$ is $B$-invariant, and $X^v$ is $B^-$-invariant). Thus, 
$X_w\cap X^v$ is equidimensional of dimension 
$\dim(X_w) + \dim(X^v) - \dim(X) = \ell(w)-\ell(v)$.
Moreover, the intersection $C_w\cap C^v$ is nonsingular and dense in
$X_w\cap X^v$. In fact, we have the following more precise result 
which may be proved by the argument of Lemma \ref{klei}; see
\cite{BrLa03} for details. 

\begin{proposition}\label{richardson}
For any $v,w\in W$, the intersection $X_w\cap X^v$ is non-empty
if and only if $v \leq w$; then $X_w\cap X^v$ is a variety.
\end{proposition}

\begin{definition}
Given $v$, $w$ in $W$ such that $v\le w$, the corresponding 
{\it Richardson variety} is $X_w^v: = X_w\cap X^v$.
\end{definition}

Note that $X_w^v$ is $T$-invariant with fixed points being the
coordinate flags $F_x = xB/B$, where $x\in W$ satisfies $v \le x \le w$. 
It follows that $X_w^v \subseteq X_{w'}^{v'}$ if and only if 
$v' \le v \le w \le w'$. Thus, the Richardson varieties may be viewed
as geometric analogues of intervals for the Bruhat order.

\begin{examples}\label{rich}
1) As special cases of Richardson varieties, we have the Schubert
varieties $X_w = X_w^{\id}$ and the opposite Schubert
varieties $X^v = X_{w_o}^v$. Also, note that the Richardson
variety $X_w^w$ is just the $T$-fixed point $F_w$, the transversal
intersection of $X_w$ and $X^w$.

\medskip

\noindent
2) Let $X_w^v$ be a Richardson variety of dimension $1$, that is,
$v\le w$ and $\ell(v) = \ell(w) -1$. Then $X_w^v$ is isomorphic to the 
projective line, and $v=ws$ for some transposition $s=s_{ij}$
(exchanging $i$ and $j$, and fixing all the other indices). More
generally, any $T$-invariant curve $Y\subset X$ is isomorphic to
$\bP^1$ and contains exactly two $T$-fixed points $v$, $w$, where
$v=ws$ for some transposition $s$. 
(Indeed, after multiplication by an element of $W$, we may assume that
$Y$ contains the standard flag $F$. Then $Y\cap C^{\id}$ is a
$T$-invariant neighborhood of $F$ in $Y$, and is also a $T$-invariant
curve in $C^{\id}\cong U^-$ (where $T$ acts by conjugation). Now any
such curve is a ``coordinate line'' given by $a_{i,j}=0$ for all
$(i,j)\ne (i_0,j_0)$, for some $(i_0, j_0)$ such that 
$1 \le j_0 < i_0 \le n$. The closure of this line in $X$ has fixed
points $F$ and $s_{i_0,j_0}F$.)
\end{examples}

Richardson varieties may be used to describe the local structure of
Schubert varieties along Schubert subvarieties, as follows. 

\begin{proposition}\label{local}
Let $v$, $w\in W$ such that $v\le w$. Then $X_w\cap vC^{\id}$ is an
open $T$-invariant neighborhood of the point $F_v$ in $X_w$, which
meets $X_w^v$ along $X_w\cap C^v$. Further, the map 
$$
(U\cap vU^-v^{-1})\times (X_w\cap C^v) \to X_w, \quad
(g,x)\mapsto gx
$$
is an open immersion with image $X_w\cap vC^{\id}$. (Recall that 
$U\cap vU^-v^{-1}$ is isomorphic to $\bC^{\ell(v)}$ as a variety, 
and that the map $U\cap vU^-v^{-1} \to X$, $g\mapsto gF_v$ 
is an isomorphism onto $C_v$.)

If, in addition, $\ell(v) = \ell(w)-1$, then $X_w\cap C^v$ is
isomorphic to the affine line. As a consequence, $X_w$ is nonsingular
along its {\rm Schubert divisor} $X_v$. 
\end{proposition}

\begin{proof}
Note that $vC^{\id}$ is an open $T$-invariant neighborhood of $F_v$ in 
$X$, isomorphic to the variety $vU^-v^{-1}$. In turn, the latter is
isomorphic to $(U\cap vU^-v^{-1})\times (U^- \cap vU^-v^{-1})$
via the product map; and the map $U^- \cap vU^-v^{-1} \to X$,
$g\mapsto gF_v$ is a locally closed immersion with image $C^v$.
It follows that the map 
$$
(U\cap vU^-v^{-1})\times C^v \to X, \quad (g,x)\mapsto gx
$$
is an open immersion with image $vF^{\id}$, and that 
$vF^{\id}\cap X^v = C^v$. Intersecting with the subvariety $X_w$ 
(invariant under the subgroup $U\cap vU^-v^{-1}$) completes the
proof of the first assertion. The second assertion follows from the
preceding example.
\end{proof}

Richardson varieties also appear when multiplying Schubert classes. 
Indeed, by Proposition \ref{richardson}, we have in $H^*(X)$:
$$
[X_w]\cup [X^v] = [X_w^v].
$$ 
Since $\dim(X_w^v) = \ell(w) - \ell(v)$, it follows that the
Poincar\'e duality pairing $\langle [X_w],[X^v]\rangle$ equals $1$ if
$w=v$, and $0$ otherwise. This implies easily the following result. 

\begin{proposition}\label{dualco}
(i) The bases $\left\{ [X_w] \right\}$ and 
$\left\{ [X^w] \right\}=\left\{[X_{w_o w}]\right\}$ 
of $H^*(X)$ are dual for the Poincar\'e duality pairing.

\noindent
(ii) For any subvariety $Y\subseteq X$, we have 
$$
[Y] = \sum_{w \in W} a^w(Y)\; [X_w],
$$ 
where 
$a^w(Y) = \langle [Y],[X^w]\rangle = \#(Y\cap gX^w)$ 
for general $g\in G$. In particular, the coefficients of $[Y]$ in the
basis of Schubert classes are non-negative.

\noindent
(iii) Let 
$$
[X_v] \cup [X_w] = \sum_{x\in W} a_{vw}^x \;
[X_x] \quad \text{in} \quad H^*(X),
$$
then the structure constants $a_{vw}^x$ are non-negative integers.
\end{proposition}

Note finally that all these results adapt readily to any partial flag
variety $G/P$. In fact, the map $f:G/B \to G/P$ induces a
ring homomorphism $f^*:H^*(G/P) \to H^*(G/B)$ which sends any Schubert
class $[X_{wP}]$ to the Schubert class $[X_{ww_{0,P}}]$, where 
$w\in W^P$. In particular, $f^*$ is injective.

\subsection{The Picard group}

In this subsection, we study the Picard group of the full flag variety
$X=G/B$. We first give a very simple presentation of this group,
viewed as the group of divisors modulo linear equivalence. The Picard
group and divisor class group of Schubert varieties will be described
in Subsection 2.2.

\begin{proposition}\label{classes}
The group $\Pic(X)$ is freely generated by the classes of the Schubert
divisors $X_{w_o s_i}$ where $i=1,\ldots,n-1$. Any ample 
(resp.~generated by its global sections) divisor on $X$ is linearly
equivalent to a positive (resp.~non-negative) combination of these
divisors. Further, any ample divisor is very ample.
\end{proposition}

\begin{proof}
The open Schubert cell $C_{w_o}$ has complement the union of the
Schubert divisors. Since $C_{w_o}$ is isomorphic to an affine
space, its Picard group is trivial. Thus, the classes of
$X_{w_o s_1},\ldots,X_{w_o s_{n-1}}$ generate the group
$\Pic(X)$. 

If we have a relation $\sum_{i=1}^{n-1} a_i X_{w_o s_i} =0$ in
$\Pic(X)$, then there exists a rational function $f$ on $X$ having a
zero or pole of order $a_i$ along each $X_{w_o s_i}$, and no other
zero or pole. In particular, $f$ is a regular, nowhere vanishing
function on the affine space $C_{w_o}$. Hence $f$ is constant, and
$a_i=0$ for all $i$.

Each Schubert divisor $X_{w_o s_d}$ is the pull-back under the
projection $X\to \Gr(d,n)$ of the unique Schubert divisor in
$\Gr(d,n)$. Since the latter divisor is a hyperplane section in the
Pl\"ucker embedding, it follows that $X_{w_o s_d}$ is generated by its
global sections. As a consequence, any non-negative combination of
Schubert divisors is generated by its global sections. Further, the
divisor $\sum_{d=1}^{n-1} X_{w_o s_d}$ is very ample, as the product
map $X \to \prod_{d=1}^{n-1} \Gr(d,n)$ is a closed immersion. 
Thus, any positive combination of Schubert divisors is very ample.

Conversely, let $D = \sum_{i=1}^{n-1} a_i X_{w_o s_i}$ be a
globally generated (resp.~ample) divisor on $X$. Then for any curve
$Y$ on $X$, the intersection number $\langle [D],[Y] \rangle$ is
non-negative (resp.~positive). Now take for $Y$ a Schubert curve
$X_{s_j}$, then
$$
\langle [D],[Y] \rangle = \langle \sum_{i=1}^{n-1}
a_i [X_{w_o s_i}], [X_{s_j}] \rangle =  
\sum_{i=1}^{n-1} a_i \langle [X^{s_i}],[X_{s_j}] \rangle = a_j.
$$
This completes the proof.
\end{proof}

\begin{remark}
We may assign to each divisor $D$ on $X$, its cohomology class 
$[D]\in H^2(X)$. Since linearly equivalent divisors are
homologically equivalent, this defines the cycle map 
$\Pic(X)\to H^2(X)$, which is an isomorphism by
Proposition \ref{classes}.

More generally, assigning to each subvariety of $X$ its cohomology
class yields the {\it cycle map} $A^*(X)\to H^{2*}(X)$, where $A^*(X)$
denotes the {\it Chow ring} of rational equivalence classes of
algebraic cycles on $X$ (graded by the codimension; in particular,
$A^1(X)= \Pic(X)$). Since $X$ has a ``cellular decomposition'' by
Schubert cells, {\it the cycle map is a ring isomorphism} by
\cite{Fu98} Example 19.1.11.

We will see in Section 4 that {\it the ring $H^*(X)$ is generated by 
$H^2(X)\cong \Pic(X)$, over the rationals}. (In fact, this holds over
the integers for the variety of complete flags, as follows easily from 
its structure of iterated projective space bundle.)
\end{remark}

Next we obtain an alternative description of $\Pic(X)$ in terms of
homogeneous line bundles on $X$; these can be defined as follows.
Let $\lambda$ be a {\it character} of $B$, i.e., a homomorphism of
algebraic groups $B \to \bC^*$. Let $B$ act on the product
$G\times\bC$ by $b (g,t): =(gb^{-1},\lambda(b)t)$. This action is
free, and the quotient 
$$
L_\lambda = G\times^B \bC :=(G\times\bC)/B
$$ 
maps to $G/B$ via $(g,t)B\mapsto gB$. This makes $L_\lambda$ the total
space of a line bundle over $G/B$, the 
{\it homogeneous line bundle associated to the weight $\lambda$}. 

Note that $G$ acts on $L_\lambda$ via $g(h,t)B:=(gh,t)B$, and that the
projection $f:L_\lambda\to G/B$ is $G$-equivariant; further, any 
$g\in G$ induces a linear map from the fiber $f^{-1}(x)$ to 
$f^{-1}(gx)$. In other words, $L_\lambda$ is a 
{\it $G$-linearized line bundle} on $X$.

We now describe the characters of $B$. Note that any such character
$\lambda$ is uniquely determined by its restriction to $T$ (since
$B=TU$, and $U$ is isomorphic to an affine space, so that any regular
invertible function on $U$ is constant). Further, one easily sees that
the characters of the group $T$ of diagonal invertible matrices are
precisely the maps
$$
\diag(t_1,\ldots,t_n)\mapsto t_1^{\lambda_1}\cdots t_n^{\lambda_n},
$$
where $\lambda_1,\ldots,\lambda_n$ are integers. This identifies the
multiplicative group of characters of $B$ (also called {\it weights})
to the additive group $\bZ^n$.

Next we express the {\it Chern classes} 
$c_1(L_\lambda)\in H^2(X)\cong \Pic(X)$ in the basis of Schubert
divisors. More generally, we obtain the {\it Chevalley formula} which
decomposes the products $c_1(L_\lambda) \cup [X_w]$ in this basis. 

\begin{proposition}\label{weights}
For any weight $\lambda$ and any $w\in W$, we have 
$$
c_1(L_\lambda) \cup [X_w] = \sum (\lambda_i-\lambda_j)\; [X_{ws_{ij}}],
$$
where the sum is over the pairs $(i,j)$ such that $1\le i<j\le n$, 
$ws_{ij}<w$, and $\ell(ws_{ij}) = \ell(w) -1$ (that is, $X_{ws_{ij}}$
is a Schubert divisor in $X_w$). In particular,
$$
c_1(L_\lambda) = \sum_{i=1}^{n-1} (\lambda_i-\lambda_{i+1}) \;
[X_{w_o s_i}] = \sum_{i=1}^{n-1} (\lambda_i-\lambda_{i+1}) \;
[X^{s_i}].
$$
Thus, the map $\bZ^n \to \Pic(X)$, $\lambda\mapsto c_1(L_\lambda)$
is a surjective group homomorphism, and its kernel is generated by
$(1,\ldots,1)$.
\end{proposition}

\begin{proof}
We may write
$$
c_1(L_\lambda) \cup [X_w] = \sum_{v\in W} a_v \; [X_v],
$$
where the coefficients $a_v$ are given by
$$
a_v = \langle c_1(L_\lambda) \cup [X_w], [X^v]\rangle
= \langle c_1(L_\lambda), [X_w]\cup [X^v]\rangle
= \langle c_1(L_\lambda), [X_w^v]\rangle.
$$
Thus, $a_v$ is the degree of the restriction of $L_\lambda$ to $X_w^v$ if 
$\dim(X_w^v)=1$, and is $0$ otherwise. Now $\dim(X_w^v)=1$ if and only
if : $v<w$ and $\ell(v)=\ell(w)-1$. Then $v=ws_{ij}$ for some
transposition $s_{ij}$, and $X_w^v$ is isomorphic to $\bP^1$, by
Example \ref{rich}.2. Further, one checks that the restriction of
$L_\lambda$ to $X_w^v$ is isomorphic to the line bundle
$\cO_{\bP^1}(\lambda_i-\lambda_j)$ of degree $\lambda_i-\lambda_j$.
\end{proof}

This relation between weights and line bundles motivates the following 

\begin{definition}
We say that the weight $\lambda =(\lambda_1,\ldots,\lambda_n)$
is {\it dominant} (resp.~{\it regular dominant}), if 
$\lambda_1 \ge \cdots \ge \lambda_n$ 
(resp.~$\lambda_1 > \cdots > \lambda_n$). 

The {\it fundamental weights} are the weights 
$\chi_1,\ldots,\chi_{n-1}$ such that 
$$
\chi_j:= (1,\ldots, 1~(j~\text{times}),
0,\ldots,0~(n-j~\text{times})).
$$
The {\it determinant} is the weight $\chi_n:=(1,\ldots,1)$.
We put 
$$
\rho := \chi_1 + \cdots + \chi_{n-1} = (n-1,n-2,\ldots,1,0).
$$
\end{definition}

By Propositions \ref{classes} and \ref{weights}, the line bundle
$L_\lambda$ is globally generated (resp.~ample) if and only if the
weight $\lambda$ is dominant (resp.~regular dominant). Further, the
dominant weights are the combinations 
$a_1 \chi_1 + \cdots + a_{n-1} \chi_{n-1} + a_n \chi_n$,
where $a_1,\ldots,a_{n-1}$ are non-negative integers, and $a_n$ is an
arbitrary integer; $\chi_n$ is the restriction to $T$ of the
determinant function on $G$. For $1\le d\le n-1$, the line bundle
$L(\chi_d)$ is the pull-back of $\cO(1)$ under the composition 
$X \to \Gr(d,n) \to \bP(\bigwedge^d \bC^n)$. Further, we have by
Proposition \ref{weights}:
$$
c_1(L_{\chi_d})\cup [X_w] = [X_{w_o s_d}] \cup [X_w] = \sum_v [X_v],
$$ 
the sum over the $v\in W$ such that $v\le w$, $\ell(v)=\ell(w)-1$, and
$v=ws_{ij}$ with $i<d<j$.

We now consider the spaces of global sections of homogeneous line
bundles. For any weight $\lambda$, we put
$$ 
H^0(\lambda): = H^0(X,L_\lambda).
$$
This is a finite-dimensional vector space, as $X$ is projective. 
Further, since the line bundle $L_\lambda$ is $G$-linearized, the
space $H^0(\lambda)$ is a {\it rational $G$-module}, i.e., $G$ acts
linearly on this space and the corresponding homomorphism 
$G \to \GL(H^0(\lambda))$ is algebraic. Further properties of this
space and a refinement of Proposition \ref{weights} are given by the
following:

\begin{proposition}\label{rep}
The space $H^0(\lambda)$ is non-zero if and only if $\lambda$ is
dominant. Then $H^0(\lambda)$ contains a unique line of eigenvectors
of the subgroup $B^-$, and the corresponding character of $B^-$ is
$-\lambda$. The divisor of any such eigenvector $p_{\lambda}$ satisfies 
$$
\div(p_{\lambda}) = \sum_{i=1}^{n-1} 
(\lambda_i - \lambda_{i+1}) \; X^{s_i}.
$$

More generally, for any $w\in W$, the $G$-module $H^0(\lambda)$
contains a unique line of eigenvectors of the subgroup $wB^-w^{-1}$,
and the corresponding weight is $-w\lambda$. Any such eigenvector
$p_{w\lambda}$ has a non-zero restriction to $X_w$, with divisor
$$
\div(p_{w\lambda}\vert_{X_w}) = 
\sum (\lambda_i - \lambda_j) \; X_{ws_{ij}},
$$ 
the sum over the pairs $(i,j)$ such that $1\le i<j\le n$ and 
$X_{ws_{ij}}$ is a Schubert divisor in $X_w$. (This makes sense as
$X_w$ is nonsingular in codimension $1$, see Proposition \ref{local}.) 

In particular, taking $\lambda=\rho$, the zero locus of
$p_{w\rho}\vert_{X_w}$ is exactly the union of all the Schubert
divisors in $X_w$. 
\end{proposition}

\begin{proof}
If $\lambda$ is dominant, then we know that $L_\lambda$ is generated
by its global sections, and hence admits a non-zero section. 
Conversely, if $H^0(\lambda)\ne 0$ then $L_\lambda$ has a section
$\sigma$ which does not vanish at some point of $X$. Since $X$ is
homogeneous, the $G$-translates of $\sigma$ generate $L_\lambda$. Thus,
$L_\lambda$ is dominant. 

Now choose a dominant weight $\lambda$ and put 
$D:=\sum_{i=1}^{n-1} (\lambda_i - \lambda_{i+1}) \; X^{s_i}$.
By Proposition \ref{weights}, we have $L_\lambda\cong \cO_X(D)$, so
that $L_\lambda$ admits a section $\sigma$ with divisor $D$. Since $D$
is $B^-$-invariant, $\sigma$ is a $B^-$-eigenvector; in particular, a
$T$-eigenvector. And since $D$ does not contain the standard flag $F$, it
follows that $\sigma(F)\ne 0$. Further, $T$ acts on the fiber of
$L_\lambda$ at $F$ by the weight $\lambda$, so that $\sigma$ has weight
$-\lambda$. If $\sigma'$ is another $B$-eigenvector in $H^0(\lambda)$,
then the quotient $\sigma'/\sigma$ is a rational function on $X$,
which is $U^-$-invariant as $\sigma$ and $\sigma'$ are. Since the
orbit $U^- F$ is open in $X$, it follows that the function
$\sigma'/\sigma$ is constant, i.e., $\sigma'$ is a scalar multiple of
$\sigma$.

By $G$-equivariance, it follows that $H^0(\lambda)$ contains a unique
line of eigenvectors of the subgroup $wB^-w^{-1}$, with weight
$-w\lambda$. Let $p_{w\lambda}$ be such an eigenvector, then
$p_{w\lambda}$ does not vanish at $F_w$, hence (by $T$-equivariance)
it has no zero on $C_w$. So the zero locus of the restriction
$p_{w\lambda}\vert_{X_w}$ has support in $X_w \setminus C_w$ and hence
is $B$-invariant. The desired formula follows by the above argument
together with Proposition \ref{weights}.
\end{proof}

\begin{remark}\label{reps}
For any dominant weight $\lambda$, the $G$-module $H^0(\lambda)$
contains a unique line of eigenvectors for $B=w_o B^- w_o$, of weight
$-w_o\lambda$. On the other hand, the evaluation of sections at
the base point $B/B$ yields a non-zero linear map 
$H^0(\lambda)\to \bC$ which is a $B$-eigenvector of weight
$\lambda$. In other words, the dual $G$-module 
$$
V(\lambda):= H^0(\lambda)^*
$$
contains a canonical $B$-eigenvector of weight $\lambda$. 

One can show that both $G$-modules $H^0(\lambda)$ and $V(\lambda)$ are
{\it simple}, i.e., they admit no non-trivial proper submodules. Further,
any simple rational $G$-module $V$ is isomorphic to $V(\lambda)$ for a
unique dominant weight $\lambda$, the {\it highest weight} of $V$. 
The $T$-module $V(\lambda)$ is the sum of its weight subspaces, and
the corresponding weights lie in the convex hull of the orbit
$W\lambda\subset \bZ^n \subset \bR^n$. For these results, see
e.g. \cite{Fu97} 8.2 and 9.3. 
\end{remark}

\begin{example}
For $d=1,\dots,n-1$, the space $\bigwedge^d \bC^n$ has a basis
consisting of the vectors
$$
e_I:=e_{i_1}\wedge\cdots\wedge e_{i_d},
$$
where $I=(i_1,\ldots,i_d)$ and $1\le i_1<\cdots<i_d\le n$. These
vectors are $T$-eigenvectors with pairwise distinct weights, and they
form a unique orbit of $W$. It follows easily that the $G$-module
$\bigwedge^d \bC^n$ is simple with highest weight $\chi_d$ (the
weight of the unique $B$-eigenvector $e_{1\ldots d}$). In other
words, we have $V(\chi_d) = \bigwedge^d \bC^n$, so that 
$H^0(\chi_d) = (\bigwedge^d \bC^n)^*$.
 
Denote by $p_I\in (\bigwedge^d \bC^n)^*$ the elements of
the dual basis of the basis $\{e_I\}$ of $\bigwedge^d \bC^n$. The
$p_I$ are homogeneous coordinates on $\Gr(d,n)$, the 
{\it Pl\"ucker coordinates}. From the previous remark, one readily
obtains that
$$
\div(p_I\vert_{X_I}) = 
\sum_{J,\; J<I,\; \vert J \vert = \vert I \vert -1} X_J.
$$
This is a refined version of the formula
$$
c_1(L) \cup [X_I] = 
\sum_{J,\;J<I,\; \vert J \vert = \vert I \vert -1} [X_J]
$$
in $H^*(\Gr(d,n))$, where $L$ denotes the pull-back of $\cO(1)$ via
the Pl\"ucker embedding. Note that $c_1(L)$ is the class of the unique
Schubert divisor.
\end{example}

\medskip

\noindent
{\bf Notes.}
The results of this section are classical; they may be found in more
detail in \cite{Fu97}, \cite{FuPr98} and \cite{Man98}, see also
\cite{Ta04}. We refer to \cite{Sp98a} Chapter 8 for an exposition of
the theory of reductive algebraic groups with some fundamental
results on their Schubert varieties. Further references are the survey
\cite{Sp98b} of Schubert varieties and their generalizations in this
setting, and the book \cite{Ku02} regarding the general framework of
Kac-Moody groups.

The irreducibility of the intersections $X_w\cap X^v$ is due to
Richardson \cite{Ri92}, whereas the intersections $C_w\cap C^v$
have been studied by Deodhar \cite{Deo85}. In fact, the Richardson
varieties in thc Grassmannians had appeared much earlier, in Hodge's 
geometric proof \cite{Ho42} of the Pieri formula which decomposes
the product of an arbitrary Schubert class with
the class of a ``special'' Schubert variety (consisting of those 
subspaces having a nontrivial intersection with a given standard 
coordinate subspace). The Richardson varieties play an important role 
in several recent articles, in relation to standard monomial theory; 
see \cite{LiSe03}, \cite{LaLi03}, \cite{KrLa03}, \cite{BrLa03}. 

The decomposition of the products $c_1(L_\lambda) \cup [X_w]$ in the
basis of Schubert classes is due to Monk \cite{Mo59} for the variety of
complete flags, and to Chevalley \cite{Ch58} in general. The Chevalley
formula is equivalent to the decomposition into Schubert classes of
the products of classes of Schubert divisors with arbitrary Schubert
classes. This yields a closed formula for certain structure constants
$a_{vw}^x$ of $H^*(X)$; specifically, those where $v=w_os_d$ for some
elementary transposition $s_d$. 

More generally, closed formulae for all the structure constants have
been obtained by several mathematicians, see \cite{KoKu86}, \cite{Du03}, 
\cite{Pr04}. The latter paper presents a general formula and applies it
to give an algebro-combinatorial proof of the Pieri formula. 
We refer to \cite{HiBo86}, \cite{PrRa93}, \cite{PrRa96}, \cite{PrRa03} 
for generalizations of the Pieri formula to the isotropic Grassmannians 
which yield combinatorial (in particular, positive) expressions for 
certain structure constants.

However, the only known proof of the positivity of the general structure 
constants is geometric. In fact, an important problem in Schubert calculus 
is to find a combinatorial expression of these constants which makes their 
positivity evident.

\vfill\eject

\section{Singularities of Schubert varieties}

As seen in Examples \ref{gras} and \ref{flag}, Schubert varieties are
generally singular. In this section, we show that their singularities
are rather mild. We begin by showing that they are normal. Then we
introduce the Bott-Samelson desingularizations, and we establish the
rationality of singularities of Schubert varieties. In particular,
these are Cohen-Macaulay; we also describe their dualizing sheaf,
Picard group, and divisor class group. Finally, we obtain the
vanishing of all higher cohomology groups $H^j(X_w,L_\lambda)$, 
where $\lambda$ is any dominant weight, and the surjectivity of the
restriction map 
$H^0(\lambda) = H^0(X,L_\lambda) \to H^0(X_w,L_\lambda)$.

\subsection{Normality}

First we review an inductive construction of Schubert cells and
varieties. Given $w\in W$ and an elementary transposition $s_i$, we
have either $\ell(s_iw)=\ell(w)-1$ (and then $s_iw < w$),
or $\ell(s_iw)=\ell(w)+1$ (and then $s_iw > w$).
In the first case, we have $Bs_i C_w=C_w \cup C_{s_iw}$,
whereas $Bs_iC_w = C_{s_iw}$ in the second case.
Further, if $w \ne \id$ (resp.~$w \ne w_o$), then there 
exists an index $i$ such that the first (resp.~second) case occurs.
(These properties of the Bruhat decomposition are easily checked in
the case of the general linear group; for arbitrary reductive groups,
see e.g. \cite{Sp98a}.)

Next let $P_i$ be the subgroup of $G=\GL_n(\bC)$ generated by $B$ and
$s_i$. (This is a {\it minimal parabolic subgroup} of $G$.) Then $P_i$ is
the stabilizer of the partial flag consisting of all the standard
coordinate subspaces, except $\langle e_1,\ldots,e_i\rangle$.
Further, $P_i/B$ is the Schubert curve $X_{s_i}\cong \bP^1$, 
and $P_i = B\cup Bs_i B$ is the closure in $G$ of $Bs_i B$.

The group $B$ acts on the product $P_i\times X_w$ by  
$b (g,x):=(gb^{-1},bx)$. This action is free; denote the quotient
by $P_i\times^B X_w$. Then the map 
$$
P_i\times X_w \to P_i\times X, \quad (g,x) \mapsto (g,gx)
$$
yields a map
$$
\iota:P_i\times^B X_w \to P_i/B\times X, \quad (g,x)B \mapsto (gB,gx).
$$
Clearly, $\iota$ is injective and its image consists of those pairs
$(gB,x)\in P_i/B\times X$ such that $g^{-1}x \in X_w$ ; this
defines a closed subset of $P_i/B\times X$. It follows that
$P_i\times^B X_w$ is a projective variety equipped with a proper
morphism 
$$
\pi:P_i\times^B X_w \to X
$$
with image $P_i X_w$, and with a morphism
$$
f:P_i\times^B X_w \to P_i/B \cong \bP^1.
$$
The action of $P_i$ by left multiplication 
on itself yields an action on $P_i\times^B X_w$; the maps
$\pi$ and $f$ are $P_i$-equivariant. Further, $f$ is a locally trivial
fibration with fiber $B\times^B X_w\cong X_w$.
 
In particular, $P_i X_w$ is closed in $X$, and hence is the closure of
$Bs_i C_w$. If $s_iw<w$, then $P_i X_w = X_w$. Then one checks
that $P_i\times^B X_w$ identifies to $P_i/B\times X_w$, so
that $\pi$ becomes the second projection. On the other hand, if 
$s_iw>w$, then $P_i X_w = X_{s_iw}$. Then 
one checks that $\pi$ restricts to an isomorphism
$$
Bs_i B\times^B C_w \to Bs_i C_w = C_{s_iw},
$$
so that $\pi$ is birational onto its image $X_{s_iw}$.

We are now in a position to prove

\begin{theorem}\label{normal}
Any Schubert variety $X_w$ is normal.
\end{theorem}

\begin{proof}
We argue by decreasing induction on $\dim(X_w) = \ell(w) =: \ell$. In
the case where $\ell=\dim(X)$, the variety $X_w = X$ is
nonsingular and hence normal. So we may assume that $\ell<\dim(X)$
and that all Schubert varieties of dimension $>\ell$ are normal.
Then we may choose an elementary transposition $s_i$ such that 
$s_iw > w$. We divide the argument into three steps.

\noindent
{\it Step 1}. We show that the morphism 
$\pi: P_i\times^B X_w\to X_{s_iw}$ 
satisfies $R^j\pi_*\cO_{P_i\times^B X_w}=0$ for all $j\ge 1$.

Indeed, $\pi$ factors as the closed immersion 
$$
\iota:P_i\times^B X_w \to 
P_i/B\times X_{s_iw} \cong \bP^1\times X_{s_iw}, \quad 
(g,x)B \mapsto (gB,gx)
$$
followed by the projection 
$$
p: \bP^1\times X_{s_iw} \to X_{s_iw}, \quad (z,x)\mapsto x.
$$
Thus, the fibers of $\pi$ are closed subschemes of $\bP^1$ and it
follows that $R^j\pi_*\cO_{P/B\times X_w}=0$ for $j>1=\dim \bP^1$.

It remains to check the vanishing of 
$R^1\pi_*\cO_{P_i\times^B X_w}$. For this, we consider the
following short exact sequence of sheaves:
$$
0 \to \cI \to \cO_{\bP^1\times X_{s_iw}} \to
\iota_*\cO_{P_i\times^B X_w}\to 0,
$$
where $\cI$ denotes the ideal sheaf of the subvariety
$P_i\times^B X_w$ of $\bP^1\times X_{s_iw}$. The derived long 
exact sequence for $p$ yields an exact sequence
$$
R^1 p_*\cO_{\bP^1\times X_{s_iw}}
\to R^1p_*(\iota_*\cO_{P_i\times^B X_w}) \to R^2 p_*\cI .
$$
Further, $R^1 p_*\cO_{\bP^1\times X_{s_iw}} =0$ as 
$H^1(\bP^1,\cO_{\bP^1})=0$;
$R^1p_*(\iota_*\cO_{P_i\times^B X_w}) =
R^1\pi_*\cO_{P_i\times^B X_w}$ as $\iota$ is a closed immersion;  
and $R^2 p_*\cI=0$ as all the fibers of $p$ have dimension $1$. This
yields the desired vanishing.

\noindent
{\it Step 2}. We now analyze the normalization map
$$
\nu:\tilde X_w\to X_w.
$$
We have an exact sequence of sheaves
$$
0 \to \cO_{X_w} \to \nu_*\cO_{\tilde X_w}\to
\cF\to 0,
$$
where $\cF$ is a coherent sheaf with support the non-normal locus of
$X_w$. Further, the action of $B$ on $X_w$ lifts to an
action on $\tilde X_w$ so that $\nu$ is equivariant. Thus, both
sheaves $\cO_{X_w}$ and $\nu_*\cO_{\tilde X_w}$ are
$B$-linearized; hence $\cF$ is $B$-linearized as well. (See
\cite{Br03b} \S 2 for details on linearized sheaves.)

Now any $B$-linearized coherent sheaf $\cG$ on $X_w$ yields an
``induced '' $P_i$-linearized sheaf $P_i\times^B \cG$ on 
$P_i\times^B X_w$ (namely, the unique $P_i$-linearized sheaf
which pulls back to the $B$-linearized sheaf $\cG$ under the inclusion  
$X_w \cong B\times^B X_w \to P_i\times^B X_w$). Further, the
assignment $\cG \mapsto P_i\times^B \cG$ is exact. 
Therefore, one obtains a short exact sequence of $P_i$-linearized
sheaves on $P_i\times^B X_w$:
$$
0 \to \cO_{P_i\times^B X_w}\to
(P_i\times^B\nu)_*\cO_{P_i\times^B\tilde X_w} \to 
P_i\times^B\cF \to 0.
$$
Apply $\pi_*$, we obtain an exact sequence of sheaves on $X_{s_iw}$:
$$
0 \to \pi_*\cO_{P_i\times^B X_w} \to
\pi_*(P_i\times^B\nu)_*\cO_{P_i\times^B\tilde X_w} \to
\pi_*(P_i\times^B\cF) \to
R^1\pi_*\cO_{P_i\times^B X_w}.
$$
Now $\pi_*\cO_{P_i\times^B X_w}=\cO_{X_{s_iw}}$ by Zariski's
main theorem, since $\pi:P_i\times^B X_w \to X_{s_iw}$ is a
proper birational morphism, and $X_{s_iw}$ is normal by the
induction assumption. Likewise, 
$\pi_*(P_i\times^B \nu)_*\cO_{P_i\times^B\tilde X_w} =
\cO_{X_{s_iw}}$. Further, $R^1\pi_*\cO_{P_i\times^B X_w}=0$
by Step 1. It follows that $\pi_*(P_i\times^B\cF)=0$.

\noindent
{\it Step 3}. Finally, we assume that $X_w$ is non-normal and we 
derive a contradiction.  

Recall that the support of $\cF$ is the non-normal locus of
$X_w$. By assumption, this is a non-empty $B$-invariant closed subset
of $X$. Thus, the irreducible components of $\supp(\cF)$ are certain
Schubert varieties $X_v$. Choose such a $v$ and let $\cF_v$
denote the subsheaf of $\cF$ consisting of sections killed by the
ideal sheaf of $X_v$ in $X_w$. Then $\supp(\cF_v)=X_v$,
since $X_v$ is an irreducible component of $\supp(\cF)$. Further, 
$\pi_*(P_i\times^B\cF_v)=0$, since $\cF_v$ is a subsheaf of
$\cF$.  

Now choose the elementary transposition $s_i$ such that
$v < s_iv$. Then $w < s_iw$ 
(otherwise, $P_iX_w=X_w$, so that $P_i$ stabilizes the 
non-normal locus of $X_w$; in particular, $P_i$ stabilizes
$X_v$, whence $s_iv<v$). Thus, the morphism
$\pi: P_i\times^B X_v \to X_{s_iv}$ restricts to an isomorphism above
$C_{s_iv}$. Since 
$\supp(P_i\times^B\cF_v) = P_i\times^B X_v$,
it follows that the support of $\pi_*(P_i\times^B\cF_v)$
contains $C_{s_iv}$, i.e., this support is the whole
$X_{s_iv}$. In particular, $\pi_*(P_i\times^B\cF_v)$ is non-zero,
which yields the desired contradiction.
\end{proof}

\subsection{Rationality of singularities}

Let $w\in W$. If $w \ne \id$ then there exists a simple
transposition $s_{i_1}$ such that $\ell(s_{i_1}w) = \ell(w)-1$.
Applying this to $s_{i_1}w$ and iterating this process, we obtain a
decomposition 
$$
w=s_{i_1}s_{i_2}\cdots s_{i_\ell}, \quad \text{where} \quad
\ell=\ell(w). 
$$
We then say that the sequence of simple transpositions 
$$
\uw:=(s_{i_1},s_{i_2},\dots,s_{i_\ell})
$$ 
is a {\it reduced decomposition} of $w$.

For such a decomposition, we have 
$X_w = P_{i_1} X_{s_{i_1}w} = P_{i_1}P_{i_2}\cdots P_{i_\ell}/B$.
We put $v:=s_{i_1}w$ and $\uv:=(s_{i_2},\dots,s_{i_\ell})$, so 
that $\uw=(s_{i_1},\uv)$ and $X_w = P_{i_1}X_v$. We define
inductively the {\it Bott-Samelson variety} $Z_{\uw}$ by 
$$
Z_{\uw} := P_{i_1}\times^B Z_{\uv}.
$$
Thus, $Z_{\uw}$ is equipped with an equivariant fibration to 
$P_{i_1}/B \cong \bP^1$ with fiber $Z_{\uv}$ at the base point.
Further, $Z_{\uw}$ is the quotient of the product 
$P_{i_1}\times\cdots\times P_{i_\ell}$ by the action of $B^\ell$
via 
$$
(b_1,\ldots,b_\ell) (g_1,g_2,\dots,g_\ell)=
(g_1 b_1^{-1},b_1 g_2 b_2^{-1},\dots,b_{\ell-1} g_\ell b_\ell^{-1}).
$$
The following statement is easily checked.

\begin{proposition}\label{bsdh}
(i) The space $Z_{\uw}$ is a nonsingular projective
$B$-variety of dimension $\ell$, where $B$ acts via 
$g (g_1,\ldots,g_\ell)B^\ell := (gg_1,\ldots,g_\ell)B^\ell$.
For any subsequence $\uv$ of $\uw$, we have a closed $B$-equivariant
immersion $Z_{\uv} \to Z_{\uw}$. 

\noindent
(ii) The map
$$
Z_{\uw}\to (G/B)^\ell =X^\ell, \quad (g_1,g_2,\dots,g_\ell)B^\ell
\mapsto (g_1B,g_1g_2B,\dots,g_1\cdots g_\ell B)
$$
is a closed $B$-equivariant embedding.

\noindent
(iii) The map 
$$
\varphi: Z_{\uw}= Z_{s_{i_1},\ldots,s_{i_\ell}} \to 
Z_{s_{i_1},\ldots,s_{i_{\ell-1}}}, \quad (g_1,\ldots,g_\ell)B^\ell
\mapsto (g_1,\ldots,g_{\ell -1})B^{\ell -1}
$$
is a $B$-equivariant locally trivial fibration with fiber
$P_{i_\ell}/B \cong \bP^1$.

\noindent
(iv) The map 
$$\pi = \pi_{\uw}: Z_{\uw}\to P_{i_1}\cdots P_{i_\ell}/B = X_w, 
\quad (g_1,\dots,g_\ell)B^\ell\mapsto g_1\cdots g_\ell B,
$$
is a proper $B$-equivariant morphism, and restricts to an isomorphism
over $C_w$. In particular, $\pi$ is birational.
\end{proposition}

An interesting combinatorial consequence of this proposition is
the following description of the Bruhat order (which may also be
proved directly).

\begin{corollary}\label{bru} 
Let $v,w\in W$. Then $v \le w$ if and only if there exist a reduced
decomposition $\uw=(s_{i_1},\dots,s_{i_\ell})$, and a subsequence
$\uv = (s_{j_1},\dots,s_{j_m})$ with product $v$. Then there exists 
a reduced subsequence $\uv$ with product $v$. 

As a consequence,  $v<w$ if and only if there exists a sequence 
$(v_1,\ldots,v_k)$ in $W$ such that $v=v_1< \cdots < v_k = w$,
and $\ell(v_{j+1}) = \ell(v_j) +1$ for all $j$.
\end{corollary}

\begin{proof}
Since $\pi$ is a proper $T$-equivariant morphism, any fiber at a
$T$-fixed point contains a fixed point (by Borel's fixed point
theorem, see e.g. \cite{Sp98a} Theorem 6.2.6). But the fixed points in 
$X_w$ (resp.~$Z_{\uw}$) correspond to the $v\in W$ such that
$v\le w$ (resp.~to the subsequences of $\uw$). This proves the first
assertion.

If $v = s_{j_1}\cdots s_{j_m}$, then the product 
$Bs_{j_1}B \cdots Bs_{j_m}B/B$ is open in $X_v$. By induction on $m$,
it follows that there exists a reduced subsequence 
$(s_{k_1},\dots,s_{k_n})$ of $(s_{j_1},\ldots,s_{j_m})$ such that
$Bs_{k_1}B \cdots Bs_{k_n}B/B$ is open in $X_v$; then 
$v = s_{k_1}\cdots s_{k_n}$. This proves the second assertion.

The final assertion follows from the second one. Alternatively,
one may observe that the complement $X_w \setminus C_w$ has pure
codimension one in $X_w$, since $C_w$ is an affine open subset of
$X_w$. Thus, for any $v<w$ there exists $x\in W$ such that 
$v\le x < w$ and $\ell(x) = \ell(w) -1$. Now induction on
$\ell(w) - \ell(v)$ completes the proof.
\end{proof}

\begin{theorem}\label{rational}
The morphism $\pi:Z_{\uw}\to X_w$ satisfies
$\pi_*\cO_{Z_{\uw}} = \cO_{X_w}$, and 
$R^j\pi_*\cO_{Z_{\uw}}=0$ for all $j\ge 1$.
\end{theorem}

\begin{proof}
We argue by induction on $\ell=\ell(w)$, the case where $\ell=0$
being trivial. For $\ell\ge 1$, we may factor $\pi=\pi_{\uw}$ as
$$
P_{i_1}\times^B \pi_{\uv} : P_{i_1}\times^B Z_{\uv} \to 
P_{i_1}\times^B X_v, \quad (g,z)B \mapsto (g,\pi_{\uv}(z))B
$$
followed by the map
$$
\pi_1: P_{i_1}\times^B X_v \to X_w, \quad (g,x) B \mapsto gx.
$$
By the induction assumption, the morphism $\pi_{\uv}$ satisfies the
conclusions of the theorem. It follows easily that so does the induced
morphism $P_{i_1}\times^B \pi_{\uv}$. But the same holds for the
morphism $\pi_1$, by the first step in the proof of Theorem
\ref{normal}. Now the Grothendieck spectral sequence for the
composition $\pi_1 \circ (P_{i_1}\times^B \pi_{\uv})=\pi_{\uw}$
(see \cite{Gro57} Chapter II) yields the desired statements.
\end{proof}

Thus, $\pi$ is a desingularization of the Schubert variety $X_w$,
and the latter has rational singularities in the following sense 
(see \cite{KeKnMuSa73} p.~49). 

\begin{definition} 
A {\it desingularization} of an algebraic variety $Y$ consists of a
nonsingular algebraic variety $Z$ together with a proper birational
morphism $\pi:Z\to Y$. We say that $Y$ has {\it rational singularities}, 
if there exists a desingularization $\pi: Z \to Y$ satisfying
$\pi_*\cO_Z=\cO_Y$ and $R^j\pi_*\cO_Z=0$ for all $j\geq 1$.
\end{definition}

Note that the equality $\pi_*\cO_Z = \cO_Y$ is equivalent to the
normality of $Y$, by Zariski's main theorem. Also, one can show that
$Y$ has rational singularities if and only if $\pi_*\cO_Z=\cO_Y$ and
$R^j\pi_*\cO_Z=0$ for all $j\geq 1$, where $\pi:Z \to Y$ is {\it any}
desigularization.

Next we recall the definition of the {\it canonical sheaf} $\omega_Y$
of a normal variety $Y$. Let $\iota:Y^\reg \to Y$ denote the inclusion
of the nonsingular locus, then $\omega_Y:=\iota_*\omega_{Y^\reg}$, where
$\omega_{Y^\reg}$ denotes the sheaf of differential forms of maximal
degree on the nonsingular variety $Y^\reg$. Since the sheaf
$\omega_{Y^\reg}$ is invertible and $\codim(Y-Y^\reg)\geq 2$, the
canonical sheaf is the sheaf of local sections of a Weil divisor
$K_Y$: the {\it canonical divisor}, defined up to linear equivalence. 
If, in addition, $Y$ is Cohen-Macaulay, then $\omega_Y$ is its 
{\it dualizing sheaf}.

For any desingularization $\pi:Z\to Y$ where $Y$ is normal, we have an 
injective trace map $\pi_*\omega_Z \to \omega_Y$. Further,
$R^j\pi_*\omega_Z=0$ for any $j\ge 1$, by the Grauert-Riemenschneider
theorem (see \cite{EnVi92} p.~59). We may now formulate the following
characterization of rational singularities, proved e.g. in
\cite{KeKnMuSa73} p.~50. 

\begin{proposition}\label{CM}
Let $Y$ be a normal variety. Then $Y$ has rational singularities if
and only if: $Y$ is Cohen-Macaulay and $\pi_*\omega_Z=\omega_Y$ for
any desingularization $\pi:Z\to Y$.
\end{proposition}

In particular, any Schubert variety $X_w$ is Cohen-Macaulay, and its
dualizing sheaf may be determined from that of a Bott-Samelson
desingularization $Z_{\uw}$. To describe the latter, put $Z:=Z_{\uw}$
and for $1\le j\le \ell$, let $Z^j\subset Z$ be the Bott-Samelson
subvariety associated with the subsequence 
$\uw^j:=(s_{i_1}\ldots,\widehat{s_{i_j}},\ldots,s_{i_\ell})$
obtained by suppressing $s_{i_j}$.

\begin{proposition}\label{clabs}
(i) With the preceding notation, $Z^1,\ldots,Z^\ell$ identify to
nonsingular irreducible divisors in $Z$, which meet transversally at a 
unique point (the class of $B^\ell$). 

\noindent
(ii) The complement in $Z$ of the boundary 
$$
\partial Z := Z^1 \cup \cdots \cup Z^\ell
$$
equals $\pi^{-1}(C_w)\cong C_w$. 

\noindent
(iii) The classes $[Z^j]$, $j=1,\ldots, \ell$, form a basis of the
Picard group of $Z$. 
\end{proposition}

Indeed, (i) follows readily from the construction of $Z$; (ii) is a
consequence of Proposition \ref{bsdh}, and (iii) is checked by the
argument of Proposition \ref{classes}. 

Next put
$$
\partial X_w := X_w \setminus C_w =
\bigcup_{v\in W, \; v<w} X_v,
$$
this is the {\it boundary} of $X_w$. By Corollary \ref{bru},
$\partial X_w$ is the union of all the Schubert divisors in $X_w$.
Further, $\pi^{-1}(\partial X_w)=\partial Z$ (as sets). 

We may now describe the dualizing sheaves of Bott-Samelson and 
Schubert varieties. 

\begin{proposition}\label{dualizing}
(i) $\omega_Z \cong (\pi^*L_{-\rho})(-\partial Z)$.

\noindent
(ii) $\omega_{X_w}\cong L_{-\rho}\vert_{X_w}(-\partial X_w)$. In
particular, $\omega_X \cong L_{-2\rho}$.

\noindent
(iii) The reduced subscheme $\partial X_w$ is Cohen-Macaulay. 
\end{proposition}

\begin{proof}
(i) Consider the curves $C_j := Z_{s_j} =P_j/B$ for $j=1,\ldots,\ell$. 
We may regard each $C_j$ as a subvariety of $Z$, namely, the transversal 
intersection of the $Z^k$ for $k\ne j$. We claim that any divisor $D$
on $Z$ such that $\langle [D], [C_j]\rangle =0$ for all $j$ is
principal. 

To see this, note that $\langle [Z^j], [C_j] \rangle =1$ for all $j$,
by Proposition \ref{clabs}. On the other hand,
$\langle [Z^j], [C_k] \rangle =0$ for all $j<k$. Indeed, we have a
natural projection 
$\varphi_j: Z = Z_{s_{i_1},\ldots,s_{i_\ell}} \to Z_{s_{i_1},\ldots,s_{i_j}}$
such that $Z^j$ is the pull-back of the corresponding divisor 
$Z_{s_{i_1},\ldots,s_{i_j}}^j$. Moreover, $\varphi_j$ maps $C_k$ to a
point whenever $k>j$. Since the $[Z^j]$ generate freely $\Pic(Z)$ by 
Proposition \ref{clabs}, our claim follows.

By this claim, it suffices to check the equality of the degrees of
the line bundles $\omega_Z(\partial Z)$ and $\pi^*L_{-\rho}$ when  
restricted to each curve $C_j$. Now we obtain 
$$
\omega_Z(\partial Z)\vert_{C_j}\cong \omega_{C_j}(\partial C_j),
$$
by the adjunction formula.
Further, $C_j\cong \bP^1$, and $\partial C_j$ is one point, so that
$\omega_{C_j}(\partial C_j)\cong \cO_{\bP^1}(-1)$. On the other hand,
$\pi$ maps $C_j$ isomorphically to the Schubert curve $X_{s_j}$, 
and $\cL_\rho\vert_{X_{s_j}} \cong \cO_{\bP^1}(1)$, so that  
$\pi^*L_{-\rho}\vert C_j \cong \cO_{\bP^1}(-1)$. This shows the
desired equality. 

(ii) Since $X_w$ has rational singularities, we have 
$\omega_{X_w} = \pi_*\omega_Z$. Further, the projection formula
yields $\omega_{X_w} \cong L_{-\rho}\otimes \pi_*\cO_Z(-\partial Z)$,
and $\pi_*\cO_Z(-\partial Z)\cong \cO_{X_w}(-\partial X_w)$ as
$\pi^{-1}(\partial X_w) = \partial Z$.

(iii) By (ii), the ideal sheaf of $\partial X_w$ in $X_w$ is locally
isomorphic to the dualizing sheaf $\omega_{X_w}$. Therefore, this
ideal sheaf is Cohen-Macaulay of depth $\dim(X_w)$. Now the exact 
sequence
$$
0 \to \cO_{X_w}(-\partial X_w) \to \cO_{X_w} \to 
\cO_{\partial X_w} \to  0
$$
yields that the sheaf $\cO_{\partial X_w}$ is Cohen-Macaulay of depth
$\dim(X_w) - 1 = \dim(\partial X_w)$.
\end{proof}

We also determine the Picard group $\Pic(X_w)$ and divisor class group
$\Cl(X_w)$ of any Schubert variety. These groups are related by an
injective map $\Pic(X_w) \to \Cl(X_w)$ which may fail to be surjective
(e.g., for $X_{24}\subset \Gr(4,2)$).

\begin{proposition}\label{class}
(i) The classes of the Schubert divisors in $X_w$ form a basis of the
divisor class group $\Cl(X_w)$.

\noindent
(ii) The restriction $\Pic(X) \to \Pic(X_w)$ is surjective, and its
kernel consists of the classes $L_\lambda$, where the weight $\lambda$
satisfies $\lambda_i = \lambda_{i+1}$ whenever $s_i\le w$. Further,
each globally generated (resp.~ample) line bundle on $X_w$ extends to
a globally generated (resp.~ample) line bundle on $X$.

\noindent
(iii) The map $\Pic(X_w) \to \Cl(X_w)$ sends the class of any
$L_\lambda$ to $\sum (\lambda_i - \lambda_j) \; X_{ws_{ij}}$ (the sum
over the pairs $(i,j)$ such that $1\le i<j\le n$ and $X_{ws_{ij}}$ is
a Schubert divisor in $X_w$).

\noindent
(iv) A canonical divisor for $X_w$ is $- \sum (j-i+1) \; X_{ws_{ij}}$ 
(the sum as above). In particular, a canonical divisor for the full
flag variety $X$ is $-2\sum_{i=1}^{n-1} X_{w_o s_i}$.
\end{proposition}

\begin{proof}
(i) is proved by the argument of Proposition \ref{classes}.

(ii) Let $L$ be a line bundle in $X_w$ and consider its pull-back
$\pi^*L$ under a Bott-Samelson desingularization 
$\pi : Z_{\uw} \to X_w$. By the argument of Proposition \ref{dualizing}, 
the class of $\pi^*L$ in $\Pic(Z_{\uw})$ is uniquely determined by its
intersection numbers $\langle c_1(\pi^*L), [C_j] \rangle$. Further, the
restriction $\pi:C_j \to \pi(C_j)$ is an isomorphism onto a Schubert
curve, and all the Schubert curves in $X_w$ arise in this way. Thus, 
$\langle c_1(\pi^*L), [C_j] \rangle$ equals either $0$ or  
$\langle c_1(L), [X_{s_i}] \rangle$ for some $i$ such that
$X_{s_i}\subseteq X_w$, i.e., $s_i\le w$. We may find a weight
$\lambda$ such that 
$\lambda_i -\lambda_{i+1} = \langle c_1(L), [X_{s_i}] \rangle$ for all
such indices $i$; then $\pi^*L_\lambda = \pi^*L$ in $\Pic(Z_{\uw})$,
whence $L=L_\lambda$ in $\Pic(X_w)$. 

If, in addition, $L$ is globally generated (resp.~ample), then 
$\langle c_1(L), [X_{s_i}] \rangle \ge 0$ (resp.~$>0$) for each
Schubert curve $X_{s_i}\subseteq X_w$. Thus, we may choose $\lambda$
to be dominant (resp.~regular dominant).

(iii) follows readily from Proposition \ref{rep}, and (iv) from
Proposition \ref{dualizing}.
\end{proof}

\subsection{Cohomology of line bundles}

The aim of this subsection is to prove the following 

\begin{theorem}\label{bwb}
Let $\lambda$ be a dominant weight and let $w\in W$. Then the
restriction map $H^0(\lambda)\to H^0(X_w,L_\lambda)$ is
surjective. Further, $H^j(X_w,L_\lambda)=0$ for any $j\ge 1$.
\end{theorem}

\begin{proof}
We first prove the second assertion in the case where $X_w = X$
is the full flag variety. Then $\omega_X \cong L_{-2\rho}$, so that
$\omega_X^{-1} \otimes L_\lambda \cong L_{\lambda+2\rho}$ is ample. 
Thus, the assertion follows from the {\it Kodaira vanishing theorem}:
$H^j(X,\omega_X\otimes \cL)=0$ for $j\ge 1$, where $\cL$ is any
ample line bundle on any projective nonsingular variety $X$.

Next we prove the second assertion for arbitrary $X_w$. For this,
we will apply a generalization of the Kodaira vanishing theorem to a
Bott-Samelson desingularization of $X_w$. Specifically, choose a
reduced decomposition $\uw$ and let $\pi:Z_{\uw}\to X_w$
be the corresponding morphism. Then the projection formula yields
isomorphisms 
$$
R^i\pi_*(\pi^* L_\lambda)\cong L_\lambda\otimes R^i\pi_*\cO_{Z_{\uw}}
$$
for all $i\ge 0$. Together with Theorem \ref{rational} and
the Leray spectral sequence for $\pi$, this yields isomorphisms 
$$
H^j(Z_{\uw},\pi^*L_\lambda)\cong H^j(X,L_\lambda)
$$ 
for all $j\ge 0$. 

We now recall a version of the 
{\it Kawamata-Viehweg vanishing theorem}, see \cite{EnVi92} \S 5. 
Consider a nonsingular projective variety $Z$, a line bundle $\cL$ on
$Z$, and a family $(D_1,\ldots,D_\ell)$ of nonsingular divisors on $Z$
intersecting transversally. Put $D := \sum_i \alpha_i D_i$, where
$\alpha_1,\ldots,\alpha_\ell$ are positive integers. Let $N$ be an
integer such that $N> \alpha_i$ for all $i$, and put $\cM:=\cL^N(-D)$.
Assume that some positive tensor power of the line bundle $\cM$ is
globally generated, and that the corresponding morphism to a
projective space is generically finite over its image (e.g., $\cM$ is
ample). Then $H^j(Z,\omega_Z\otimes \cL)=0$ for all $j\ge 1$.

We apply this result to the variety $Z := Z_{\uw}$, the line bundle
$\cL := (\pi^*L_{\lambda+\rho})(\partial Z)$, and the divisor
$D := \sum_i (N-b_i) Z^i$ where $b_1,\ldots,b_\ell$ are positive
integers such that $\sum_i b_i Z^i$ is ample (these exist by Lemma
\ref{ggample} below). Then 
$\cL^N(-D) = (\pi^*L_{N(\lambda+\rho)})(b_1Z^1+\cdots+b_\ell Z^\ell)$
is ample, and $\omega_Z\otimes\cL=\pi^*L_\lambda$. This yields the
second assertion.

To complete the proof, it suffices to show that the restriction map
$H^0(X_w,L_\lambda)\to H^0(X_v,L_\lambda)$ is surjective whenever 
$w = s_iv > v$ for some elementary transposition $s_i$.
As above, this reduces to checking the surjectivity of the 
restriction map 
$H^0(Z,\pi^*L_\lambda) \to H^0(Z^1,\pi^*L_\lambda)$. 
For this, by the long exact sequence
$$ 
0\to H^0(Z,(\pi^*L_\lambda)(-Z^1)) \to H^0(Z,\pi^*L_\lambda) 
\to H^0(Z^1,\pi^*L_\lambda) \to H^1(Z,(\pi^*L_\lambda)(-Z^1)),
$$
it suffices in turn to show the vanishing of 
$H^1(Z,(\pi^*L_\lambda)(-Z^1))$. 

We will deduce this again from the Kawamata-Viehweg vanishing theorem.
Let $a_1$, $\ldots$, $a_\ell$ be positive integers such that the line
bundle $(\pi^*L_{\lambda+a_1\rho})(a_2 Z^2 + \cdots + a_\ell Z^\ell)$ is
ample (again, these exist by Lemma \ref{ggample} below).
Put $\cL:= (\pi^*L_{\lambda+\rho})(Z^2 + \cdots+ Z^\ell)$ and 
$D:=\sum_{i=2}^\ell (N-a_i)Z^i$, where $N>a_1,a_2,\ldots,a_\ell$. Then
$\cL^N(-D)= (\pi^*L_{N(\lambda+\rho)})(a_2 Z^2+\cdots+ a_\ell Z^\ell)$ 
is ample, and $\omega_Z\otimes\cL = (\pi^*L_\lambda)(-Z^1)$. Thus, we
obtain $H^j(Z,(\pi^*L_\lambda)(-Z^1))=0$ for all $j\ge 1$.
\end{proof}

\begin{lemma}\label{ggample}
Let $Z = Z_{\uw}$ with boundary divisors $Z^1,\ldots,Z^{\ell}$. Then
there exist positive integers $a_1,\ldots,a_\ell$ such that the line
bundle $(\pi^*L_{a_1\rho})(a_2 Z^2 + \cdots + a_\ell Z^\ell)$ is ample. 
Further, there exist positive integers $b_1,\ldots,b_\ell$ such that 
the divisor $b_1 Z^1 + \cdots + b_\ell Z^\ell$ is ample.
\end{lemma}

\begin{proof}
We prove the first assertion by induction on $\ell$. If $\ell=1$, then 
$\pi$ embeds $Z$ into $X$, so that $\pi^*L_{a_1\rho}$ is ample for any
$a_1>0$. In the general case, the map 
$$
\varphi:Z\to Z^\ell = 
(P_{i_1}\times \cdots \times P_{i_{\ell-1}})/B^{\ell -1}, 
\quad  
(g_1,\ldots,g_\ell)B^\ell\mapsto (g_1,\ldots,g_{\ell-1})B^{\ell-1}
$$
fits into a cartesian square
$$\CD
Z @>{\varphi}>> Z^\ell\\
@V{\pi}VV @V{\psi}VV\\
G/B @>{f}>> G/P_{i_\ell},\\
\endCD$$
where $\psi((g_1,\ldots,g_{\ell-1})B^{\ell-1}) = 
g_1\cdots g_{\ell-1}P_{i_\ell}$.
Further, the boundary divisors $Z^{1,\ell}$, $\ldots$, $Z^{\ell-1,\ell}$
of $Z^\ell$ satisfy $\varphi^*Z^{i,\ell}=Z^i$. Denote by 
$$
\pi_\ell : Z^\ell = Z_{(s_{i_1},\ldots,s_{i_{\ell-1}})} 
\to X_{s_{i_1}\cdots s_{i_{\ell-1}}} = X_{w s_{i_\ell}}
$$ 
the natural map. By the induction assumption, there exist positive
integers $a_1,a_2,\ldots,a_{\ell-1}$ such that the line bundle
$(\pi_\ell^*L_{a_1\rho})(a_2 Z^{1,\ell} + \cdots + 
a_{\ell-1} Z^{\ell-1,\ell})$
is very ample on $Z^\ell$. Hence its pull-back
$$
\varphi^*((\pi_\ell^*L_{a_1\rho})(a_2 Z^{2,\ell} + \cdots + 
a_{\ell-1} Z^{\ell-1,\ell}))
= (\varphi^*\pi_\ell^*L_{a_1\rho})(a_2 Z^2 + \cdots + 
a_{\ell-1} Z^{\ell-1})
$$
is a globally generated line bundle on $Z$. Thus, it suffices to show
that the line bundle 
$\pi^*L_{b\rho}\otimes (\varphi^*\pi_\ell^*L_{-a_1\rho})(a_1Z^\ell)$
is globally generated and $\varphi$-ample for $b\gg a_1$. 
(Indeed, if $\cM$ is a globally generated, $\varphi$-ample line bundle
on $Z$, and $\cN$ is an ample line bundle on $Z^\ell$, then 
$\cM\otimes \varphi^*\cN$ is ample on $Z$). Equivalently, it suffices
to show that 
$\pi^*L_{c\rho}\otimes \pi^*L_\rho \otimes 
(\varphi^*\pi^*L_{-\rho})(Z^\ell)$ 
is globally generated and $\varphi$-ample for $c\gg 0$. But we have
by Proposition \ref{clabs}: 
$$
\pi^*L_{\rho}\otimes (\varphi^*\pi_\ell^*L_{-\rho})(Z^\ell)
=\omega_Z^{-1}(-\partial Z)\otimes 
\varphi^*(\omega_{Z^\ell}(\partial Z^\ell))(Z^\ell)
= \omega_Z^{-1}\otimes \varphi^*\omega_{Z^\ell}
= \omega_\varphi^{-1} = \pi^*\omega_f^{-1},
$$
where $\omega_\varphi$ (resp.~$\omega_f$) denotes the relative
dualizing sheaf of the morphism $\varphi$ (resp.~$f$).
Further, $L_{c\rho}\otimes\omega_f^{-1}$ is very ample on $G/B$ for 
$c\gg 0$, as $L_\rho$ is ample. Thus,
$\pi^*(L_{c\rho}\otimes \omega_f^{-1})$ is globally generated and
$\varphi$-ample. This completes the proof of the first assertion.

The second assertion follows by recalling that the restriction of
$L_\rho$ to $X_w$ admits a section vanishing exactly on 
$\partial X_w$ (Remark \ref{reps} 2). Thus, $\pi^*L_\rho$ admits a
section vanishing exactly on $\partial Z = Z^1\cup\cdots\cup Z^\ell$. 
\end{proof}

Next we consider a regular dominant weight $\lambda$ and the
corresponding very ample homogeneous line bundle $L_\lambda$.
This defines a projective embedding
$$
X \to \bP(H^0(X,L_\lambda)^*) = \bP(V(\lambda))
$$
and, in turn, a subvariety $\tilde X\subseteq V(\lambda)$, invariant
under the action of $G\times \bC^*$, where $\bC^*$ acts by scalar
multiplication. We say that $\tilde X$ is the 
{\it affine cone over $X$} associated with this projective
embedding. Likewise, we have the affine cones $\tilde X_w$ over
Schubert varieties.

\begin{corollary}\label{projnorm}
For any regular dominant weight $\lambda$, the affine cone over
$X_w$ in $V(\lambda)$ has rational singularities. In particular,
$X_w$ is projectively normal in its embedding into
$\bP(V(\lambda))$.
\end{corollary}

\begin{proof}
Consider the total space $Y_w$ of the line bundle 
$L_{\lambda}^{-1}\vert_{X_w}$. We have a proper morphism
$$
\pi:Y_w \to \tilde X_w
$$
which maps the zero section to the origin, and restricts to an
isomorphism from the complement of the zero section to the complement
of the origin. In particular, $\pi$ is birational. Further, $Y_w$ has
rational singularities, since it is locally isomorphic to 
$X_w \times \bC$. Thus, it suffices to show that 
the natural map $\cO_{\tilde X_w} \to \pi_*\cO_{Y_w}$ is surjective,
and $R^j\pi_*\cO_{Y_w}=0$ for any $j\ge 1$. Since $\tilde X_w$ is
affine, this amounts to: the algebra $H^0(Y_w,\cO_{Y_w})$ is generated
by the image of $H^0(\lambda)$, and $H^j(Y_w,\cO_{Y_w})=0$ for 
$j\ge 1$. Further, since the projection $f:Y_w \to X_w$ is affine and
satisfies
$$
f_*\cO_{Y_w} = \bigoplus_{n=0}^{\infty} L_\lambda^{\otimes n} = 
\bigoplus_{n=0}^{\infty} L_{n\lambda},$$ 
we obtain
$$
H^j(Y_w,\cO_{Y_w}) = \bigoplus_{n=0}^{\infty} H^j(X_w, L_{n\lambda}).
$$
So $H^j(Y_w,\cO_{Y_w})=0$ for $j\ge 1$, by Theorem \ref{bwb}. To
complete the proof, it suffices to show that  
the algebra $\bigoplus_{n=0}^{\infty} H^0(X_w, L_{n\lambda})$ is
generated by the image of $H^0(\lambda)$. Using the surjectivity of
the restriction maps  
$H^0(L_{n\lambda}) \to H^0(X_w, L_{n\lambda})$
(Theorem \ref{bwb} again), it is enough to consider the case where
$X_w = X$. Now the multiplication map
$$
H^0(\lambda)^{\otimes n} \to H^0(n\lambda), \quad
\sigma_1 \otimes \cdots \otimes \sigma_n \mapsto 
\sigma_1\cdots \sigma_n
$$
is a non-zero morphism of $G$-modules. Since $H^0(n\lambda)$ is
simple, this morphism is surjective, which completes the proof.
\end{proof}

\medskip

\noindent
{\bf Notes.} 
In their full generality, the results of this section were obtained
by many mathematicians during the mid-eighties. Their most elegant
proofs use reduction to positive characteristics and the techniques of 
Frobenius splitting, see \cite{MeRa85}, \cite{RaRa85}, \cite{Ra87}. 

Here we have presented alternative proofs: for normality and
rationality of singularities, we rely on an argument of Seshadri
\cite{Se87} simplified in \cite{Br03b}, which is also valid in
arbitrary characteristics. For cohomology of line bundles, our
approach (based on the Kawamata-Viehweg vanishing theorem) is a
variant of that of Kumar; see \cite{Ku02}. 

The construction of the Bott-Samelson varieties is due to $\ldots$
Bott and Samelson \cite{BoSa58} in the framework of compact Lie groups,
and to Hansen \cite{Han73} and Demazure \cite{Dem74} in our
algebro-geometric setting. The original construction of Bott and
Samelson is also presented in \cite{Du04} with applications to the
multiplication of Schubert classes.

The line bundles on Bott-Samelson varieties
have been studied by Lauritzen and Thomsen in \cite{LaTh04}; in
particular, they determined the globally generated (resp.~ample) line
bundles. On the other hand, the description of the Picard group and
divisor class group of Schubert varieties is due to Mathieu in
\cite{Mat88}; it extends readily to any Schubert variety $Y$ in any
flag variety $X = G/P$. One may also show that the boundary of $Y$ is
Cohen-Macaulay, see \cite{Br03b} Lemma 4. But a simple formula for the
dualizing sheaf of $Y$ is only known in the case where $X$ is the full
flag variety.

An important open question is the explicit determination of the
singular locus of a Schubert variety, and of the corresponding generic
singularities (i.e., the singularities along each irreducible
component of the singular locus). The book \cite{BiLa00} by Billey and
Lakshmibai is a survey of this question, which was recently solved 
(independently and simultaneously) by several mathematicians in the
case of the general linear group; see \cite{BiWa03}, \cite{Co03},
\cite{KaLaRe03}, \cite{Man01a}, \cite{Man01b}. The generic
singularities of Richardson varieties are also worth investigating.

\vfill\eject

\section{The diagonal of a flag variety}

Let $X=G/B$ be the full flag variety and denote by $\diag(X)$ the
diagonal in $X\times X$. In this section, we construct 
a degeneration of $\diag(X)$ in $X\times X$ to the union of all the
products $X_w\times X^w$, where the $X_w$ (resp.~$X^w$) are the
Schubert (resp.~opposite Schubert) varieties.

Specifically, we construct a subvariety 
$\cX\subseteq X\times X\times \bP^1$ such that the fiber of the
projection $\pi : \cX \to \bP^1$ at any $t\ne 0$ is isomorphic to
$\diag(X)$, and we show that the fiber at $0$ (resp.~$\infty$) is the
union of all the $X_w\times X^w$ (resp.~$X^w\times X_w$). For this, we
use the normality of $\cX$ which is deduced from a general
normality criterion for varieties with group actions, obtained in turn
by adapting the argument for the normality of Schubert varieties.

Then we turn to applications to the Grothendieck ring $K(X)$. After a
brief presentation of the definition and main properties of
Grothendieck rings, we obtain two additive bases of $K(X)$ which are
dual for the bilinear pairing given by the Euler characteristic of the
product. Further applications will be given in Section 4.

\subsection{A degeneration of the diagonal}

We begin by determining the cohomology class of $\diag(X)$ in 
$X\times X$, where $X$ is the full flag variety.

\begin{lemma}\label{diagco}
We have $[\diag(X)] = \sum_{w\in W} [X_w\times X^w]$ in 
$H^*(X\times X)$.
\end{lemma}

\begin{proof}
By the results in Subsection 1.3 and the K\"unneth isomorphism, a
basis for the abelian group $H^*(X\times X)$ consists of the classes 
$[X_w\times X^v]$, where $v,w\in W$. Further, the dual basis 
(with respect to the Poincar\'e duality pairing) consists of the  
$[X^w\times X_v]$. Thus, we may write
$$
[\diag(X)]=\sum_{v,w\in W} a_{wv}\;[X_w\times X^v],
$$ 
where the coefficients $a_{wv}$ are given by
$$
a_{wv} = \langle [\diag(X)], [X^w\times X_v] \rangle.
$$
Further, since $X^w$ meets $X_v$ properly along $X^w_v$ with
intersection multiplicity $1$, it follows that $\diag(X)$ meets
$X^w\times X_v$ properly along $\diag(X_v^w)$ in $X\times X$ with
intersection multiplicity $1$. This yields
$$ 
[\diag(X)] \cup [X^w \times X_v] = [\diag(X_w^v)].
$$  
And since $\dim(X_v^w)=0$ if and only if $v=w$, we see that 
$a_{wv}$ equals $1$ if $v=w$, and $0$ otherwise.
\end{proof}

This formula suggests the existence of a degeneration of $\diag(X)$ to 
$\bigcup_{w\in W} X_w\times X^w$. We now construct such a
degeneration. The idea is to move $\diag(X)$ in $X\times X$ by a
general one-parameter subgroup of the torus $T$ acting on $X\times X$
via its action on the second copy, and to take limits.  

Specifically, let 
$$
\lambda:\bC^*\to T, \quad t\mapsto(t^{a_1},\dots,t^{a_n})
$$
where $a_1,\ldots, a_n$ are integers satisfying
$a_1>\cdots>a_n$. Define $\cX$ to be the closure in  
$X\times X\times \bP^1$ of the subset 
$$
\{(x,\lambda(t)x,t) ~\vert~ x\in X, \; t\in\bC^*\}\subseteq
X\times X\times\bC^*.
$$
Then $\cX$ is a projective variety, and the fibers of the projection 
$\pi : \cX \to \bP^1$ identify with closed subschemes of $X\times X$. 
Further, the fiber $\pi^{-1}(1)$ equals $\diag(X)$. In fact,
$\pi^{-1}(\bC^*)$ identifies to $\diag(X) \times \bC^*$ via
$(x,y,t)\mapsto (x,\lambda(t^{-1})x,t)$, and this identifies 
the restriction of $\pi$ to the projection 
$\diag(X) \times \bC^* \to \bC^*$. 

\begin{theorem}\label{fibers}
We have equalities of subschemes of $X\times X$:
$$
\pi^{-1}(0)=\bigcup_{w\in W} X_w\times X^w
\quad \text{and} \quad
\pi^{-1}(\infty)=
\bigcup_{w\in W} X^w\times X_w.
$$
\end{theorem}

\begin{proof}
By symmetry, it suffices to prove the first equality.
We begin by showing the inclusion 
$\bigcup_{w\in W} X_w\times X^w\subseteq \pi^{-1}(0)$.
Equivalently, we claim that $C_w\times C^w\subset \pi^{-1}(0)$ for all
$w\in W$.

For this, we analyze the structure of $X\times X$ in a neighborhood of
the base point $(F_w,F_w)$ of $C_w\times C^w$ (recall that $F_w$
denotes the image under $w$ of the standard flag $F$). By
Proposition \ref{local}, $wC^{\id}$ is a $T$-invariant
open neighborhood of $F_w$ in $X$, isomorphic to $w U^-w^{-1}$. 
Further, $C_w = U  F_w \cong (wU^-w^{-1}\cap U) F_w$ identifies via
this isomorphism to the subgroup $w U^-w^{-1}\cap U$. Likewise, $C^w$
identifies to the subgroup $w U^-w^{-1}\cap U^-$, and the product
map in the group $wU^- w^{-1}$
$$
(wU^-w^{-1}\cap U) \times (w U^-w^{-1}\cap U^-) \to w U^-w^{-1}
$$
is an isomorphism. Further, each factor is isomorphic to an affine space.

The group $\bC^*$ acts on $w U^-w^{-1}$
via its homomorphism $t\mapsto (t^{a_1},\ldots,t^{a_n})$ to $T$ and the
action of $T$ on $w U^-w^{-1}$ by conjugation. In fact, this action
is linear, and hence $w U^-w^{-1}$ decomposes into a direct sum of weight
subspaces. Using the assumption that $a_1>\cdots>a_n$, one checks that
the sum of all the positive weight
subspaces is $wU^-w^{-1}\cap U=C_w$; likewise, the sum of all the
negative weight subspaces is $C^w$. In other words,
$$
C_w=\{x\in w U^-w^{-1} ~\vert~ \lim_{t\to 0}\lambda(t)x=\id\},
\quad
C^w=\{y\in w U^-w^{-1} ~\vert~ \lim_{t\to\infty}\lambda(t)y=\id\}.
$$
Now identify our neighborhood $wC^{\id}\times wC^{\id}$
with $C_w\times C^w\times C_w\times C^w$. Take arbitrary
$x\in C_w$ and $y\in C^w$, then
$$
(x,\lambda(t)^{-1}y,\lambda(t)x,y)\to(x,\id,\id,y) 
\quad \text{as} \quad t\to 0.
$$
By the definition of $\cX$, it follows that $\pi^{-1}(0)$
contains the point $(x,\id,\id,y)$, identified to $(x,y)\in X\times X$.
This proves the claim.

From this claim, it follows that $\pi^{-1}(0)$ contains 
$\bigcup_{w\in W} X_w\times X^w$ (as schemes). On the other hand, the
cohomology class of $\pi^{-1}(0)$ equals that of $\pi^{-1}(1)$, i.e.,
$\sum_{w\in W} [X_w\times X^w]$ by Lemma \ref{diagco}. 
Further, the cohomology class of any non-empty subvariety of 
$X\times X$ is a positive integer combination of classes 
$[X_w \times X^v]$ by Proposition \ref{dualco}. It follows that the
irreducible components of $\pi^{-1}(0)$ are exactly the 
$X_w\times X^w$, and that the corresponding multiplicities
are all $1$. Thus, the scheme $\pi^{-1}(0)$ is generically reduced.

To complete the proof, it suffices to show that $\pi^{-1}(0)$ is
reduced. Since $\pi$ may be regarded as a regular function on 
$\cX$, it suffices in turn to show that $\cX$ is normal.
In the next subsection, this will be deduced from a general normality
criterion for varieties with group action.
\end{proof}

To apply Theorem \ref{fibers}, we will also need to analyze the
structure sheaf of the special fiber $\pi^{-1}(0)$.
This is the content of the following statement. 

\begin{proposition}\label{filt}
The sheaf $\cO_{\pi^{-1}(0)}$ admits a filtration with
associated graded
$$
\bigoplus_{w\in W} \cO_{X_w}\otimes\cO_{X^w}
(-\partial X^w).
$$
\end{proposition}

\begin{proof}
We may index the finite partially ordered set $W=\{w_1,\ldots,w_N\}$
so that $i\le j$ whenever $w_i\le w_j$ (then $w_N = w_o$). Put
$$
Z_i:=X_{w_i}\times X^{w_i} \quad \text{and} \quad
Z_{\ge i}:=\bigcup_{j\ge i} Z_j,
$$
where $1\le i\le N$. Then $Z_{\ge 1} = \pi^{-1}(0)$ and 
$Z_{\ge N} = X_{w_o}\times X^{w_o} = X\times \{w_oF\}$.
Further, the $Z_{\ge i}$ form a decreasing filtration of
$\pi^{-1}(0)$. 
This yields exact sequences
$$
0\to \cI_i \to \cO_{Z_{\ge i}} \to \cO_{Z_{\ge i+1}} \to 0,
$$
where $\cI_i$ denotes the ideal sheaf of $Z_{\ge i+1}$ in $Z_{\ge i}$.
In turn, these exact sequences yield an increasing filtration of 
the sheaf $\cO_{\pi^{-1}(0)}$ with associated graded 
$\bigoplus_i \cI_i$. Since 
$Z_{\ge i} = Z_{\ge i+1} \cup Z_i$, we may identify $\cI_i$ with the
ideal sheaf of $Z_i\cap Z_{\ge i+1}$ in $Z_i= X_{w_i}\times X^{w_i}$. 
To complete the proof, it suffices to show that
$$
Z_i\cap Z_{\ge i+1} = X_{w_i}\times \partial X^{w_i}.
$$
We first check the inclusion ``$\subseteq$''. Note that 
$Z_i\cap Z_{\ge i+1}$ is invariant under $B\times B^-$, and hence is a
union of products $X_u\times X^v$ for certain $u,v\in W$. We must have
$u\le w_i\le v$ (since $X_u\times X^v\subseteq Z_i$) and $w_i\ne v$
(since $X_u\times X^v\subseteq Z_{\ge i+1}$). Thus, 
$X_u\times X^v\subseteq X_{w_i}\times \partial X^{w_i}$. 
To check the opposite inclusion, note that if 
$X^v\subseteq \partial X^{w_i}$ then $v>w_i$, so that $v=w_j$ with 
$j>i$. Thus, 
$X_{w_i}\times X^v\subset X_{w_j}\times X^{w_j}\subseteq Z_{\ge i+1}$.
\end{proof}

\subsection{A normality criterion}

Let $G$ be a connected linear algebraic group acting on an
algebraic variety $Z$. Let $Y\subset Z$ be a subvariety, 
invariant under the action of a Borel subgroup $B\subseteq G$, and let
$P\supset B$ be a parabolic subgroup of $G$. Then, as in Subsection
2.1, we may define the ``induced'' variety $P\times^B Y$. It is
equipped with a $P$-action and with $P$-equivariant maps
$\pi:P\times^B Y \to Z$ (a proper morphism with image $PY$), and
$f:P\times^B Y\to P/B$ (a locally trivial fibration with fiber $Y$). 
If, in addition, $P$ is a minimal parabolic subgroup (i.e.,
$P/B\cong\bP^1$), and if $PY\ne Y$, then $\dim(PY)=\dim(Y)+1$, and the
morphism $\pi$ is generically finite over its image $PY$.

We say that $Y$ is {\it multiplicity-free} if it satisfies 
the following conditions:

\noindent
(i) $G Y=Z$.

\noindent
(ii) Either $Y=Z$, or $Z$ contains no $G$-orbit.  

\noindent
(iii) For all minimal parabolic subgroups $P\supset B$ such that
$PY\neq Y$, the morphism $\pi:P\times^B Y\to PY$ is birational, 
and the variety $PY$ is multiplicity-free.

(This defines indeed the class of multiplicity-free subvarieties
by induction on the codimension, starting with $Z$).

For example, Schubert varieties are multiplicity-free. Further, the
proof of their normality given in Subsection 2.1 readily adapts to
show the following

\begin{theorem}\label{multfree}
Let $Y$ be a $B$-invariant subvariety of a $G$-variety $Z$. If $Z$ is
normal and $Y$ is multiplicity-free, then $Y$ is normal. 
\end{theorem}

Next we obtain a criterion for multiplicity-freeness of any $B$-stable
subvariety of $Z:=G$, where $G$ acts by left multiplication. 
Note that the $B$-stable subvarieties $Y\subseteq G$ correspond to 
the subvarieties $V$ of the full flag variety $G/B$, by taking
$V := \{g^{-1}B ~\vert~ g\in Y\}$.

\begin{lemma}\label{01}
With the preceding notation, $Y$ is multiplicity-free if and only if 
$[V]$ is a multiplicity-free combination of Schubert classes, i.e., 
the coefficients of $[V]$ in the basis $\{[X_w]\}$ are either
$0$ or $1$. Equivalently, $\langle [V],[X^w]\rangle \le 1$ for all $w$.
\end{lemma}

\begin{proof}
Clearly, $Y$ satisfies conditions (i) and (ii) of
multiplicity-freeness. For condition (iii), consider a minimal
parabolic subgroup $P\supset B$ and the natural map $f:G/B\to G/P$. 
Then the subvariety of $G$ associated with $f^{-1}f(V)$ is $PY$. As a
consequence, $PY\ne Y$ if and only if the restriction 
$f\vert_V:V\to f(V)$ is generically finite. Further, the fibers of
$f\vert_V$ identify to those of the natural map 
$\pi:P\times^B Y\to PY$; in particular, both  morphisms have the same
degree $d$. Note that $d=1$ if and only if $\pi$ (or, equivalently,
$f\vert_V$) is birational.

Let $X_w\subseteq G/B$ be a Schubert variety of positive dimension. 
We may write $X_w = P_1\cdots P_\ell/B$, where $(P_1,\ldots,P_\ell)$
is a sequence of minimal parabolic subgroups, and $\ell=\dim(X_w)$. 
Put $P:=P_\ell$ and $X_v:=P_1\cdots P_{\ell-1}/B$. Then 
$X_w= f^{-1}f(X_w)$, and the restriction $X_v\to f(X_v)=f(X_w)$ is
birational. We thus obtain the equalities of intersection numbers
$$\displaylines{
\langle [V], [X_w] \rangle_{G/B} = 
\langle [V], f^{-1}[f(X_w)] \rangle_{G/B} = 
\langle f_*[V], [f(X_w)] \rangle_{G/P} 
\hfill\cr\hfill
= d \; \langle [f(V)], [f(X_w)] \rangle_{G/P} =
d \; \langle [f^{-1}f(V), [X_v] \rangle_{G/B},
\cr}$$
as follows from the projection formula and from the equalities
$f_*[V] = d \; [f(V)]$, $f_*[X_v] = [f(X_v)] = [f(X_w)]$. 
From these equalities, it follows that $[V]$ is a multiplicity-free
combination of Schubert classes if and only if: $d=1$ and 
$[f^{-1}f(V)]$ is a multiplicity-free combination of Schubert classes,
for any minimal parabolic subgroup $P$ such that $PY\ne Y$. Now the
proof is completed by induction on $\codim_{G/B}(V)=\codim_G(Y)$.
\end{proof}

We may now complete the proof of Theorem \ref{fibers} by showing that
$\cX$ is normal. Consider first the group $G\times G$, the Borel
subgroup $B\times B$, and the variety $Z: = G\times G$, where 
$G\times G$ acts by left multiplication. Then the
subvariety $Y := (B\times B)\diag(G)$ is multiplicity-free.
(Indeed, $Y$ corresponds to the variety $V = \diag(X) \subset X\times X$,
where $X=G/B$. By Lemma \ref{diagco}, the coefficients of
$[\diag(X)]$ in the basis of Schubert classes are either $0$ or $1$,
so that Lemma \ref{01} applies.)

Next consider the same group $G\times G$ and take 
$Z := G\times G\times\bP^1$, where $G\times G$ acts via left
multiplication on the factor $G\times G$. Let $Y$ be the preimage in
$Z$ of the subvariety $\cX \subset X\times X\times \bP^1$
under the natural map 
$G\times G\times \bC \to X\times X\times \bP^1$
(a locally trivial fibration). Clearly, $Y$ satisfies conditions
(i), (ii) of multiplicity-freeness. Further, condition (iii) follows
from the fact that $Y$ contains an open subset isomorphic to 
$(B\times B)\diag(G)\times \bC^*$, together with the
multiplicity-freeness of $(B\times B)\diag(G)$. Since $Z$ is
nonsingular, it follows that $Y$ is normal by Theorem
\ref{multfree}. Hence, $\cX$ is normal as well.

\subsection{The Grothendieck group}

For any scheme $X$, the {\it Grothendieck group of coherent
sheaves on $X$} is the abelian group $K_0(X)$ generated by symbols 
$[\cF]$, where $\cF$ is a coherent sheaf on $X$, subject to the
relations $[\cF] = [\cF_1] + [\cF_2]$ whenever there exists an exact
sequence of sheaves 
$0\to \cF_1 \to \cF \to \cF_2\to 0$.
(In particular, $[\cF]$ only depends on the isomorphism class of $\cF$.)
For example, any closed subscheme $Y\subseteq X$ yields a class
$[\cO_Y]$ in $K(X)$.

Likewise, we have the {\it Grothendieck group $K^0(X)$ of vector
bundles on $X$}, generated by symbols $[E]$, where $E$ is a vector
bundle on $X$, subject to the relations $[E] = [E_1] + [E_2]$
whenever there exists an exact sequence of vector bundles 
$0\to E_1\to E \to E_2\to 0$. The tensor product of vector bundles
yields a commutative, associative multiplication law on
$K^0(X)$ denoted by $(\alpha,\beta)\mapsto \alpha \cdot \beta$. With
this multiplication, $K^0(X)$ is a commutative ring, the identity
element being the class of the trivial bundle of rank $1$.

The duality of vector bundles $E\mapsto E^{\vee}$ is compatible with
the defining relations of $K^0(X)$. Thus, it yields a map 
$K^0(X) \to K^0(X)$, $\alpha \mapsto \alpha^{\vee}$, which is an
involution of the ring $K^0(X)$: the {\it duality involution.}
 
By associating with each vector bundle $E$ its (locally free) sheaf of
sections $\cE$, we obtain a map 
$$
\varphi:K^0(X) \to K_0(X).
$$ 
More generally, since tensoring with a locally free sheaf is exact, the
ring $K^0(X)$ acts on $K_0(X)$ via 
$$
[E]\cdot [\cF] := [\cE\otimes_{\cO_X} \cF],
$$ 
where $E$ is a vector bundle on $X$ with sheaf of sections $\cE$, and
$\cF$ is a coherent sheaf on $X$. 
This makes $K_0(X)$ a module over $K^0(X)$; further, 
$\varphi(\alpha) = \alpha \cdot [\cO_X]$ for any $\alpha \in K^0(X)$.

If $Y$ is another scheme, then the external tensor product of sheaves
(resp.~vector bundles) yields product maps 
$K_0(X)\times K_0(Y) \to K_0(X\times Y)$, 
$K^0(X)\times K^0(Y) \to K^0(X\times Y)$,
compatible with the corresponding maps $\varphi$. We will denote both
product maps by $(\alpha,\beta)\mapsto \alpha\times \beta$.

If $X$ is a nonsingular variety, then $\varphi$ is an isomorphism. 
In this case, we identify $K_0(X)$ with $K^0(X)$, and
we denote this ring by $K(X)$, the {\it Grothendieck ring of $X$}.
For any coherent sheaves $\cF$, $\cG$ on $X$, we have
$$
[\cF] \cdot [\cG] = \sum_j (-1)^j \; [Tor_j^X(\cF,\cG)].
$$
(This formula makes sense because the sheaves $Tor_j^X(\cF,\cG)$
are coherent, and vanish for $j > \dim(X)$). In particular, 
$[\cF] \cdot [\cG]=0$ if the sheaves $\cF$ and $\cG$ have disjoint 
supports. Further, 
$$
[\cF]^{\vee} = \sum_j (-1)^j \; [Ext_X^j(\cF,\cO_X)].
$$
In particular, if $Y$ is an equidimensional Cohen-Macaulay subscheme
of $X$, then 
$$
[\cO_Y]^{\vee} = (-1)^c \; [Ext_X^c(\cO_Y,\cO_X)]
= (-1)^c \; [\omega_{Y/X}] 
= (-1)^c \; [\omega_Y] \cdot [\omega_X]^{\vee},
$$
where $c$ denotes the codimension of $Y$,and 
$\omega_{Y/X}: = \omega_Y \otimes \omega_X^{-1}$ denotes the 
relative dualizing sheaf of $Y$ in $X$.

Returning to an arbitrary scheme $X$, any morphism of schemes
$f:X\to Y$ yields a pull-back map 
$$
f^*:K^0(Y)\to K^0(X), \quad [E]\mapsto [f^*E].
$$
If, in addition, $f$ is {\it flat}, then it defines similarly a
pull-back map $f^*:K_0(Y)\to K_0(X)$. 

On the other hand, any {\it proper} morphism $f:X\to Y$ yields a
push-forward map
$$
f_*:K_0(X)\to K_0(Y), \quad [\cF]\mapsto \sum_j (-1)^j \; [R^j f_*(\cF)].
$$
As above, this formula makes sense as the higer direct images
$R^j f_*(\cF)$ are coherent sheaves on $Y$, which vanish for $j>\dim(X)$.
Moreover, we have the {\it projection formula}
$$
f_*((f^*\alpha) \cdot \beta) = \alpha \cdot f_*\beta
$$ 
for all $\alpha\in K^0(Y)$ and $\beta\in K_0(X)$.

In particular, if $X$ is complete then we obtain a map 
$$
\chi:K_0(X)\to \bZ, \quad 
[\cF] \mapsto \chi(\cF) = \sum_j (-1)^j \; h^j(\cF),
$$
where $h^j(\cF)$ denotes the dimension of the $j$-th cohomology
group of $\cF$, and $\chi$ stands for the {\it Euler-Poincar\'e
characteristic}.

We will repeatedly use the following result of ``homotopy invariance''
in the Grothendieck group.

\begin{lemma}\label{const}
Let $X$ be a variety and let $\cX$ be a subvariety of $X\times\bP^1$
with projections $\pi:\cX\to\bP^1$ and $p:\cX\to X$. Then the class
$[\cO_{p(\pi^{-1}(z))}]\in K_0(X)$ is independent of $z \in \bP^1$. 
\end{lemma}

\begin{proof}
The exact sequence 
$0\to \cO_{\bP^1}(-1) \to \cO_{\bP^1} \to \cO_z\to 0$ of
sheaves on $\bP^1$ shows that the class $[\cO_z]\in K_0(\bP^1)$ is
independent of $z$. Since $\pi$ is flat, it follows that the class 
$\pi^*[\cO_z] = [\cO_{\pi^{-1}(z)}]\in K_0(\cX)$
is also independent of $z$, and the same holds for 
$p_*[\cO_{\pi^{-1}(z)}] \in K_0(X)$ since $p$ is proper. But  
$p_*[\cO_{\pi^{-1}(z)}]= [\cO_{p(\pi^{-1}(z))}]$, since 
$p$ restricts to an isomorphism $\pi^{-1}(z)\to p(\pi^{-1}(z))$.
\end{proof}

Finally, we present a relation of $K_0(X)$ to the Chow group $A_*(X)$
of rational equivalence classes of algebraic cycles on $X$
(graded by the dimension), see \cite{Fu98} Example 15.1.5. Define
the {\it topological filtration} on $K_0(X)$ by letting $F_j K_0(X)$
to be the subgroup generated by coherent sheaves whose support has
dimension at most $j$. Let $\Gr K_0(X)$ be the associated graded
group. Then assigning to any subvariety $Y\subseteq X$ the class
$[\cO_Y]$ passes to rational equivalence (as
follows from Lemma \ref{const}) and hence defines a morphism
$A_*(X)\to \Gr K_0(X)$ of graded abelian groups. This morphism is
surjective ; it is an isomorphism over the rationals if, in addition,
$X$ is nonsingular (see \cite{Fu98} Example 15.2.16).

\subsection{The Grothendieck group of the flag variety}

The Chow group of the full flag variety $X$ is isomorphic to its
cohomology group and, in particular, is torsion-free. It follows that
the associated graded of the Grothendieck group (for the topological
filtration) is isomorphic to the cohomology group; this isomorphism
maps the image of the structure sheaf $\cO_Y$ of any subvariety, to
the cohomology class $[Y]$. Thus, the following result may be viewed
as a refinement in $K(X\times X)$ of the equality 
$[\diag(X)] =\sum_{w\in W} [X_w\times X^w]$ in $H^*(X\times X)$.

\begin{theorem}\label{dualk}
(i) In $K(X\times X)$ holds 
$$
[\cO_{\diag(X)}] =  \sum_{w\in W}
[\cO_{X_w}]\times [\cO_{X^w}(-\partial X^w)].
$$

\noindent
(ii) The bilinear map 
$$
K(X)\times K(X)\to \bZ, \quad
(\alpha,\beta)\mapsto \chi(\alpha \cdot \beta)
$$ 
is a nondegenerate pairing. Further, 
$\{[\cO_{X_w}]\}$, $\{[\cO_{X^w}(-\partial X^w)]\}$
are bases of the abelian group $K(X)$, dual for this pairing.
\end{theorem}

\begin{proof}
(i) By Theorem \ref{fibers} and Lemma \ref{const}, we have
$[\cO_{\diag(X)}] = [\cO_{\bigcup_{w\in W} X_w\times X^w}]$.
Further, 
$[\cO_{\bigcup_{w\in W} X_w\times X^w}]=\sum_{w\in W}
[\cO_{X_w}] \times [\cO_{X^w}(-\partial X^w)]$ 
by Proposition \ref{filt}. 

(ii) Let $p_1,p_2:X\times X\to X$ be the projections. Let $\cE$ be a
locally free sheaf on $X$. Then we have by (i):
$$\displaylines{
[\cE\vert_{\diag(X)}] = [p_2^*\cE]\cdot [\cO_{\diag(X)}]
\hfill\cr\hfill
=\sum_{w\in W} [p_2^*\cE]\cdot 
[p_1^*\cO_{X_w}\otimes p_2^*\cO_{X^w}(-\partial X^w)]
= \sum_{w\in W} [p_1^* \cO_{X_w}\otimes p_2^* \cE\vert_{X^w}(-\partial X^w)].
\cr}$$
Applying $(p_1)_*$ to both sides and using the projection formula yields
$$
[\cE]=\sum_{w\in W}\chi(\cE\vert_{X^w}(-\partial X^w)) \; [\cO_{X_w}].
$$
Since the group $K(X)$ is generated by classes of locally free
sheaves, it follows that 
$$
\alpha = \sum_{w\in W}\chi(\alpha\cdot [\cO_{X^w}(-\partial X^w)]) 
\;[\cO_{X_w}]
$$
for all $\alpha \in K(X)$. Thus, the classes $[\cO_{X_w}]$ generate the
group $K(X)$.

To complete the proof, it suffices to show that these classes are
linearly independent. If $\sum_{w\in W} n_w \; [\cO_{X_w}]=0$ is a
non-trivial relation in $K(X)$, then we may choose $v\in W$ maximal
such that $n_v\ne 0$. Now a product $[\cO_{X_w}]\cdot [\cO_{X^v}]$ is
non-zero only if $X_w\cap X^v$ is non-empty, i.e., $v\le w$. Thus, we
have by maximality of $v$: 
$$
0 = \sum_{w\in W} n_w [\cO_{X_w}]\cdot [\cO_{X^v}] =
n_v [\cO_{X_v}]\cdot [\cO_{X^v}].
$$
Further, we have $[\cO_{X_v}]\cdot [\cO_{X^v}]=[\cO_{vF}]$. (Indeed,
$X_v$ and $X^v$ meet transversally at the unique point $vF$; see
Lemma \ref{fp} below for a more general result). Further, $[\cO_{vF}]$
is non-zero since $\chi(\cO_{vF})=1$; a contradiction.
\end{proof}

We put for simplicity
$$
\cO_w:=[\cO_{X_w}] \quad \text{and} \quad
\cI_w:=[\cO_{X_w}(-\partial X_w)].
$$
The $\cO_w$ are the {\it Schubert classes} in $K(X)$. Further,
$\cI_w = [\cO_{X_w}] - [\cO_{\partial X_w}]$ by the exact sequence 
$0 \to \cO_{X_w}(-\partial X_w) \to \cO_{X_w} \to \cO_{\partial X_w}
\to 0$. We will express the $\cI_w$ in terms of the $\cO_w$, and vice
versa, in Proposition \ref{oi} below.

Define likewise
$$
\cO^w:=[\cO_{X^w}] \quad \text{and} \quad 
\cI^w:=[\cO_{X^w}(-\partial X^w)].
$$
In other words, $\cO^w = [\cO_{w_0X_{w_0w}}]$ and 
$\cI^w = [\cO_{w_0 X_{w_0w}}(-w_0\partial X_{w_0w})]$. But
$[\cO_{gY}] = [\cO_Y]$ for any $g\in G$ and any subvariety 
$Y\subseteq X$. Indeed, this follows from Lemma \ref{const} together
with the existence of a connected chain of rational curves in
$G$ joining $g$ to $\id$ (since the group $G$ is generated by images
of algebraic group homomorphisms $\bC\to G$ and $\bC^*\to G$). Thus,
$$
\cO^w = \cO_{w_0w} \quad \text{and} \quad \cI^w = \cI_{w_0w}.
$$  
Now Theorem \ref{dualk}(ii) yields the equalities  
$$
\alpha = \sum_{w\in W} \chi(\alpha \cdot \cI^w) \; \cO_w 
= \sum_{w\in W} \chi(\alpha \cdot \cO_w) \; \cI^w,
$$ 
for any $\alpha\in K(X)$. 

\begin{remarks}
1) Theorem \ref{dualk} and the isomorphism $\Gr K(X)\cong H^*(X)$
imply that the classes $\cO_w$ ($w\in W$, $\ell(w)\le j$) form a basis
of $F_j K(X)$; another basis of this group consists of the $\cI_w$ 
($w\in W$, $\ell(w)\le j$).

\medskip

\noindent
2) All the results of this section extend to an arbitrary flag variety
$G/P$ by replacing $W$ with the set $W^P$ of minimal representatives.
\end{remarks}

\begin{examples}\label{db}
1) Consider the case where $X$ is the projective space $\bP^n$. 
Then the Schubert varieties are the linear subspaces $\bP^j$, 
$0\le j\le n$, and the corresponding opposite Schubert varieties
are the $\bP^{n-j}$. Further, $\partial \bP^j = \bP^{j-1}$ so that 
$$
[\cO_{\bP^j}(-\partial \bP^j)]= [\cO_{\bP^j}] - [\cO_{\bP^{j-1}}]
=[\cO_{\bP^j}(-1)].
$$
Thus, $\{[\cO_{\bP^j}]\}$ is a basis of $K(\bP^n)$ with dual
basis $\{[\cO_{\bP^{n-j}}(-1)]\}$. 

The group $K(\bP^n)$ may be described more concretely in terms of
polynomials, as follows. For each coherent sheaf $\cF$ on $\bP^n$,
the function $\bZ\to \bZ$, $k\mapsto \chi(\cF(k))$ is
polynomial of degree equal to the dimension of the support of $\cF$;
this defines the {\it Hilbert polynomial} 
$P_{\cF}(t) \in \bQ [t]$. Clearly, 
$P_{\cF}(t) = P_{\cF_1}(t) + P_{\cF_2}(t)$ for any exact sequence
$0\to \cF_1 \to \cF \to \cF_2 \to 0$. Thus, the Hilbert polynomial
yields an additive map 
$$
P: K(\bP^{n-1})\to \bQ[t], \quad
[\cF] \mapsto P_{\cF}(t).
$$
Since $\chi(\cO_{\bP^j}(k))= {k+j\choose j}$, it follows that $P$ maps
the basis $\{[\cO_{\bP^j}]\}$ to the linearly independent polynomials
$\{ {t+j\choose j} \}$. Thus, $P$ identifies $K(\bP^{n-1})$ to the 
additive group of polynomials of degree $\le n$ in one variable which
take integral values at all integers. Note that $P$ takes non-zero
values at classes of non-trivial sheaves.

\medskip

\noindent
2) More generally, consider the case where $X$ is a Grassmannian.
Let $L$ be the ample generator of $\Pic(X)$, then the boundary of each
Schubert variety $X_I$ (regarded as a reduced Weil divisor on $X_I$) is
the divisor of the section $p_I$ of $L\vert_{X_I}$; see Remark
\ref{reps}.3. Thus, we have an exact sequence 
$$
0 \to L^{-1}\vert_{X_I} \to  \cO_{X_I} \to \cO_{\partial X_I} \to 0,
$$
where the map on the left is the multiplication by $p_I$. It follows
that
$$
[\cO_{X_I}(-\partial X_I)] = [L^{-1}\vert_{X_I}].
$$
Thus, the dual basis of the basis of Schubert classes 
$\{\cO_{X_I} := \cO_I\}$ is the basis $\{[L^{-1}]\cdot \cO^I\}$.
\end{examples}

\medskip

\noindent
{\bf Notes.}
The cohomology class of the diagonal is discussed in \cite{FuPr98}
Appendix G, in a relative situation which yields a generalization of
Lemma \ref{diagco}.

Our degeneration of the diagonal of a flag variety was first
constructed in \cite{BrPo00}, by using canonical compactifications of
adjoint semisimple groups; see \cite{Br03a} for further developments
realizing these compactifications as irreducible components of Hilbert
schemes. The direct construction of 3.1 follows \cite{BrLa03}
with some simplifications. In [loc.cit.], this degeneration was
combined with vanishing theorems for unions of Richardson varieties,
to obtain a geometric approach to standard monomial theory. 
Conversely, this theory also yields the degeneration of the diagonal
presented here, see \cite{LaLi03}.

The normality criterion in 3.2 appears first in \cite{Br03b}. It is
also proved there that a $B$-invariant multiplicity-free subvariety
$Y$ of a $G$-variety $Z$ is normal and Cohen-Macaulay (resp. has
rational singularities), if $Y$ is normal and Cohen-Macaulay
(resp. has rational singularities). This yields an alternative proof
for the rationality of singularities of Schubert varieties. 

The exposition in 3.3 is based on \cite{BoSe58} regarding
fundamental results on the Grothendieck ring $K(X)$, where $X$ is any
nonsingular variety, and on \cite{Fu98} regarding the relation of this
ring to intersection theory on $X$. The reader will find another
overview of $K$-theory in \cite{Bu04} together with several
developments concerning degeneracy loci. In particular, a combinatorial 
expression for the structure constants of the Grothendieck ring of 
Grassmannians is presented there, after \cite{Bu02}. This yields 
another proof of the result in Example \ref{db}(ii); see the proof 
of Corollary 1 in \cite{Bu04}.

The dual bases of the $K$-theory of the flag manifold presented in 3.4
appear in \cite{LaSc83} for the variety of complete flags. In the
general framework of $T$-equivariant $K$-theory of flag varieties,
they were constructed by Kostant and Kumar \cite{KoKu90}. In fact, our
approach fits into this framework. Indeed, $T$ acts on 
$X\times X\times \bP^1$ via $t (x,y,z) = (tx,ty,z)$. This action
commutes with the $\bC^*$-action via $\lambda$ and leaves $\cX$
invariant; clearly, the morphism $\pi: \cX \to \bP^1$ is $T$-invariant
as well. Thus, $\pi$ is a degeneration of $T$-varieties. Further, the
filtration of $\cO_{\pi^{-1}}(0)$ constructed in Proposition
\ref{filt} is also $T$-invariant. So Theorem \ref{dualk} extends
readily to the $T$-equivariant Grothendieck group. 

The idea of determining the (equivariant) class of a subvariety by an
(equivariant) degeneration to a union of simpler subvarieties
plays an essential role in the articles of Graham \cite{Gra01} on the
structure constants of the equivariant cohomology ring of flag
varieties, and of Knutson and Miller \cite{KnMi01} on Schubert
polynomials. These polynomials are special representatives of Schubert
classes in the cohomology ring of the variety of complete flags. They
were introduced by Lascoux and Sch\"utzenberger \cite{LaSc82},
\cite{LaSc85} and given geometric interpretations in 
\cite{Gol01}, \cite{KnMi01}. Likewise, the Grothendieck polynomials are
special representatives of Schubert classes in the Grothendieck
ring of the complete flag variety, see \cite{LaSc83} and \cite{Bu04}.

It would be very interesting to have further examples of
varieties with a torus action, where the diagonal admits an
equivariant degeneration to a reduced union of products of
subvarieties. The Bott-Samelson varieties should provide such
examples; their $T$-equivariant Grothendieck ring has been described 
by Willems \cite{Wi04} with applications to equivariant Schubert
calculus that generalize results of Duan \cite{Du03}.

\vfill\eject

\section{Positivity in the Grothendieck group of the flag variety}

Let $Y$ be a subvariety of the full flag variety $X=G/B$. By the
results of Section 3, we may write in the Grothendieck group $K(X)$:
$$
[\cO_Y] = \sum_{w\in W} c^w(Y) \; \cO_w,
$$
where the $\cO_w = [\cO_{X_w}]$ are the Schubert classes. Further,
$c^w(Y)=0$ unless $\dim(Y)\ge \dim(X_w)=\ell(w)$, and we have in the
cohomology group $H^*(X)$:
$$
[Y] = \sum_{w\in W, \; \ell(w) = \dim(Y)} c^w(Y) \; [X_w].
$$
By Proposition \ref{dualco}, it follows that 
$c^w(Y) = \#(Y\cap gX^w)$ for general $g\in G$, if 
$\ell(w) = \dim(Y)$; in particular, $c^w(Y)\ge 0$ in this case. 

One may ask for the signs of the integers $c^w(Y)$, where $w$ is
arbitrary. In this section, we show that these signs are 
{\it alternating}, i.e., 
$$
(-1)^{\dim(Y)-\ell(w)} c^w(Y) \ge 0,
$$
whenever $Y$ has rational singularities (but not for arbitrary $Y$, 
see Remark \ref{surf}.2). 

We also show that the Richardson varieties have rational singularities,
and we generalize to these varieties the results of Section 2
for cohomology groups of homogeneous line bundles on Schubert
varieties. From this, we deduce that the structure constants of the
ring $K(X)$ in its basis of Schubert classes have alternating signs as
well, and we present several related positivity results.

Finally, we obtain a version in $K(X)$ of the Chevalley formula, that
is, we decompose the product $[L_\lambda] \cdot \cO_w$ in the basis
of Schubert classes, where $\lambda$ is any dominant weight, and $X_w$
is any Schubert variety.

\subsection{The class of a subvariety}

In this subsection, we sketch a proof of the alternation of signs for
the coefficients $c^w(Y)$. By Theorem \ref{dualk}, we have
$$
c^w(Y) = \chi([\cO_Y]\cdot [\cO_{X^w}(-\partial X^w)])
= \chi([\cO_Y]\cdot [\cO_{X^w}]) - 
\chi([\cO_Y]\cdot [\cO_{\partial X^w}]).
$$
Our first aim is to obtain a more tractable formula for $c^w(Y)$.
For this, we need the following version of a lemma of Fulton and
Pragacz (see \cite{FuPr98} p.~108).

\begin{lemma}\label{fp}
Let $Y$, $Z$ be equidimensional Cohen-Macaulay subschemes of a
nonsingular variety $X$. If $Y$ meets $Z$ properly in $X$, then   
the scheme-theoretic intersection $Y\cap Z$ is equidimensional and
Cohen-Macaulay, of dimension $\dim(Y)+\dim(Z)-\dim(X)$. Further,
$$
Tor_i^X(\cO_Y,\cO_Z) = 0 = Tor_i^X(\omega_Y,\omega_Z)
$$
for any $j \ge 1$, and 
$\omega_{Y\cap Z}= \omega_Y\otimes \omega_Z\otimes \omega_X^{-1}$.

Thus, we have in $K(X)$:
$$
[\cO_{Y\cap Z}] = [\cO_Y] \cdot [\cO_Z] 
\quad \text{and} \quad [\omega_{Y\cap Z}] = 
[\omega_Y] \cdot [\omega_Z] \cdot[\omega_X^{-1}].
$$
\end{lemma}

We also need another variant of Kleiman's transversality theorem
(Lemma \ref{klei}):

\begin{lemma}\label{rs}
Let $Y$ be a Cohen-Macaulay subscheme of the flag variety $X$ and let
$w\in W$. Then $Y$ meets properly $gX^w$ for general $g\in G$; further,
$Y\cap gX^w$ is equidimensional and Cohen-Macaulay. 

If, in addition, $Y$ is a variety with rational singularities, then
$Y\cap gX^w$ is a disjoint union of varieties with rational
singularities (again, for general $g\in G$).
\end{lemma}

We refer to \cite{Br03b} p. 142--144 for the proof of these results. 
Together with the fact that the boundary of any Schubert variety is
Cohen-Macaulay (Corollary \ref{dualizing}), they imply that
$$\displaylines{
c^w(Y) = \chi(\cO_{Y\cap gX^w}) - \chi(\cO_{Y\cap g\partial X^w})
\hfill\cr\hfill
=\chi(\cO_{Y\cap gX^w}(-Y\cap g\partial X^w))
= \sum_{j=0}^{\dim(Y\cap gX^w)} 
(-1)^j \; h^j(\cO_{Y\cap gX^w}(-Y\cap g\partial X^w)).
\cr}$$
Further, 
$\dim(Y\cap gX^w) = \dim(Y) + \dim(X^w) -\dim(X) = \dim(Y) - \ell(w)$.
Thus, the assertion on the sign of $c^w(Y)$ will result from the
following vanishing theorem, which holds in fact for any partial flag variety
$X$.

\begin{theorem}\label{va}
Let $Y\subseteq X$ be a subvariety with rational singularities and let
$Z \subseteq X$ be a Schubert variety. Then we have for general 
$g\in G$: 
$$
H^j(Y\cap gZ, \cO_{Y\cap gZ}(-Y\cap g\partial Z))=0 \quad
\text{whenever} \quad j<\dim(Y) + \dim(Z) -\dim(X).
$$
\end{theorem}

\begin{proof}
First we present the argument in the simplest case, where $X=\bP^n$
and $Y$ is nonsingular. Then $Z=\bP^j$ and 
$\cO_Z(-\partial Z)= \cO_{\bP^j}(-1)$, see Example \ref{db}.1. 
Thus, $Y\cap gZ=:V$ is a general linear section of $Y$. By Bertini's
theorem, $V$ is nonsingular (and irreducible if its dimension is
positive). Further,
$\cO_{Y\cap gZ}(-Y\cap g\partial Z)=\cO_V(-1)$. Thus, we are reduced 
to showing the vanishing of $H^j(V,\cO(-1))$ for $j<\dim(V)$, where $V$
is a nonsingular subvariety of $\bP^n$. But this follows from the
Kodaira vanishing theorem.

Next we consider the case where $X$ is a Grassmannian, and $Y$ is
allowed to have rational singularities. Let $L$ be the ample
generator of $\Pic(X)$ and recall that 
$\cO_Z(-\partial Z) = L^{-1}\vert_Z$. It follows that 
$\cO_{Y\cap gZ}(-Y\cap g\partial Z)= L^{-1}\vert_{Y\cap gZ}$. 
Further, by Lemma \ref{rs}, $Y\cap gZ$ is a disjoint union of 
varieties with rational singularities, of dimension
$\dim(Y)+\dim(Z)-\dim(X)$. Thus, it suffices to show that
$H^j(V,L^{-1})=0$ whenever $V$ is a variety with rational
singularities, $L$ is an ample line bundle on $V$, and $j<\dim(V)$. 
Let $\pi:\tilde V\to V$ be a desingularization and put 
$\tilde L:=\pi^*L$. Since $R^i \pi_*\cO_{\tilde V}=0$ for any 
$i\ge 1$, we obtain isomorphisms 
$H^j(V,L^{-1})\cong H^j(\tilde V, \tilde L^{-1})$ for all $j$. 
Thus, the Grauert-Riemenschneider theorem (see \cite{EnVi92} Corollary
5.6) yields the desired vanishing. 

The proof for arbitrary flag varieties goes along similar lines, but
is much more technical. Like in the proof of Theorem \ref{bwb},
one applies the Kawamata-Viehweg theorem to a desingularization of
$Y\cap gZ$; see \cite{Br03b} p.~153--156 for details.
\end{proof}

\begin{remarks}\label{surf}
1) As a consequence of Theorem \ref{va}, we have 
$$
c^w(Y) = (-1)^{\dim(Y) - \ell(w)} \; 
h^{\dim(Y) - \ell(w)}(\cO_{Y\cap gX^w}(-Y\cap g\partial X^w)).
$$
By using Serre duality on $Y\cap gX^w$, it follows that
$$
c^w(Y) = (-1)^{\dim(Y) - \ell(w)} \; 
h^0(Y\cap gX^w, L_\rho \otimes \omega_Y).
$$

\medskip

\noindent
2) The property of alternation of signs for the coefficients of
$[\cO_Y]$ on Schubert varieties fails for certain (highly singular)
subvarieties $Y$ of a flag variety $X$. Indeed, there exist surfaces 
$Y\subset X=\bP^4$ such that the coefficient of $[\cO_Y]$ on $[\cO_x]$
(where $x$ is any point of $\bP^4$) is arbitrarily negative.

Specifically, let $d\ge 3$ be an integer and let $C$ be the image of
the morphism 
$\bP^1\to \bP^3$, $(x,y)\mapsto (x^d,x^{d-1}y,xy^{d-1},y^d)$
(a closed immersion). Then $C$ is a nonsingular rational curve of
degree $d$ in $\bP^3$. Regarding $C$ as a curve in
$\bP^4\supset\bP^3$, choose $x \in \bP^4 \setminus \bP^3$ and denote by
$Y\subset \bP^4$ the projective cone over $C$ with vertex $x$, that
is, the union of all projective lines containing $x$ and meeting $C$.
Then $Y$ is a surface, so that we have by Example \ref{db}.1:
$$
[\cO_Y] = c_2(Y) \; [\cO_{\bP^2}] + c_1(Y) \; [\cO_{\bP^1}] + c_0(Y)
\; [\cO_x].
$$
We claim that $c_0(Y) \le 3 - d$.

To see this, first notice that $c_0(Y) = \chi(\cO_Y(-1))$, as 
$\chi(\cO_{\bP^j}(-1))=0$ for all $j\ge 1$. Thus, 
$$
c_0(Y) =\chi(\cO_Y)-\chi(\cO_{Y\cap \bP^3})
=\chi(\cO_Y)-\chi(\cO_C)=\chi(\cO_Y)-1.
$$
To compute $\chi(\cO_Y)$, consider the desingularization $\pi:Z\to Y$,
where $Z$ is the total space of the projective line bundle
$\bP(\cO_C \oplus \cO_C(-1))$ on $C$ (that is, the blow-up of $x$ in $Y$). 
Then we have an exact sequence
$$
0\to \cO_Y \to \pi_*\cO_Z \to \cF \to 0,
$$
where the sheaf $\cF$ is supported at $x$. Further, $R^i\pi_*\cO_Z=0$
for all $i\ge 1$. (Indeed, since the affine cone $Y_0:=Y\setminus C$
is an affine neighborhood of $x$ in $Y$, it suffices to show that
$H^i(Z_0,\cO_{Z_0})=0$ for $i\ge 1$, where $Z_0:=\pi^{-1}(Y_0)$. Now 
$Z_0$ is the total space of the line bundle 
$\cO_C(-1)\cong \cO_{\bP^1}(-d)$ on $C\cong\bP^1$, whence 
$$
H^i(Z_0,\cO_{Z_0})\cong \bigoplus_{n=0}^{\infty}
H^i(\bP^1,\cO_{\bP^1}(nd))
$$
for any $i\ge 0$.) Thus, we obtain 
$\chi(\cO_Y) = \chi(\cO_Z) - \chi(\cF) = 1- h^0(\cF)$, so that 
$c_0(Y) = -h^0(\cF)$. Further, $\cF$ identifies with the quotient
$(\pi_*\cO_{Z_0})/\cO_{Y_0}$. Since $Y_0\subset \bC^4$ is the affine
cone over $C\subset \bP^3$, this quotient is a graded vector space
with component of degree $1$ being 
$H^0(\cO_{\bP^1}(d))/H^0(\cO_{\bP^3}(1))$, of dimension $d-3$. Thus, 
$h^0(\cF)\ge d-3$. This completes the proof of the claim.

On the other hand, for any surface $Y\subset \bP^n$, the coefficient
$c_2(Y)$ is the degree of $Y$, a positive integer. 
Further, one checks that 
$$
c_1(Y) = \chi(\cO_Y(-1)) - \chi(\cO_Y(-2)) = 
\chi(\cO_{Y\cap \bP^{n-1}}(-1)) = -h^1(\cO_{Y\cap \bP^{n-1}}(-1))
$$
for any hyperplane $\bP^{n-1}$ which does not contain $Y$. Thus,
$c_1(Y)\le 0$.

Likewise, one may check that the property of alternation of signs
holds for any curve in any flag variety. In other words, the
preceding counterexample has the smallest dimension.
\end{remarks}

\subsection{More on Richardson varieties}

We begin with a vanishing theorem for these varieties that
generalizes Theorem \ref{bwb}. Let $v$, $w$ in $W$ such that $v\le w$
and let $X_w^v$ be the corresponding Richardson variety. Then $X_w^v$
has two kinds of boundaries, namely
$$
(\partial X_w)^v := (\partial X_w)\cap X^v \quad \text{and} \quad
(\partial X^v)_w := (\partial X^v)\cap X_w,
$$
where $\partial X^v = X^v \setminus C^v = \bigcup_{u>v} X^u$ denotes
the boundary of the opposite Schubert variety $X^v$. Define the total
boundary by 
$$
\partial X_w^v:= (\partial X_w)^v \cup (\partial X^v)_w,
$$
this is a closed subset of pure codimension $1$ in $X_w^v$. We may now 
state

\begin{theorem}\label{ri}
(i) The Richardson variety $X_w^v$ has rational singularities, and its
dualizing sheaf equals $\cO_{X_w^v}(-\partial X_w^v)$. Further, we
have in $K(X)$: 
$$
[\cO_{X_w^v}] = \cO_w \cdot \cO^v = \cO_w \cdot \cO_{w_o v}.
$$

\noindent
(ii) $H^j(X_w^v,L_\lambda)=0$ for any $j\ge 1$ and any dominant weight
$\lambda$.

\noindent
(iii) $H^j(X_w^v,L_\lambda(-(\partial X^v)_w))=0$ 
for any $j\ge 1$ and any dominant weight $\lambda$.

\noindent
(iv) $H^j(X_w^v,L_\lambda(-\partial X_w^v))=0$ 
for any $j\ge 1$ and any regular dominant weight $\lambda$.
\end{theorem}

\begin{proof}
(i) follows from the rationality of singularities of Schubert
varieties and the structure of their dualizing sheaves, together with
Lemmas \ref{fp} and \ref{rs}.

(ii) We adapt the proof of Theorem \ref{bwb} to this setting. 
Choose a reduced decomposition $\uw$ of $w$ and let $Z_{\uw}$ be the
associated Bott-Samelson variety with morphism
$$
\pi_{\uw}: Z_{\uw}\to X_w.
$$
Likewise, a reduced decomposition $\uv$ of $v$ yields an
opposite Bott-Samelson variety $Z^{\uv}$ (defined via the opposite
Borel subgroup $B^-$ and the corresponding minimal parabolic
subgroups) together with a morphism
$$
\pi^{\uv}: Z^{\uv}\to X^v.
$$
Now consider the fibered product
$$
Z = Z_{\uw}^{\uv} := Z_{\uw} \times_X Z^{\uv}
$$
with projection 
$\pi = \pi_{\uw}^{\uv}: Z_{\uw}^{\uv}\to X_w \cap X^v = X_w^v$.
Using Kleiman's transversality theorem, one checks that
$Z_{\uw}^{\uv}$ is a nonsingular variety and $\pi$ is
a desingularization of $X_w^v$. Let $\partial Z$ be the union of the
boundaries
$$
(\partial Z_{\uw})^{\uv} := \partial Z_{\uw} \times_X Z^{\uv},
\quad
(\partial Z^{\uv})_{\uw} := Z_{\uw} \times_X \partial Z^{\uv}.
$$ 
This is a union of irreducible nonsingular divisors intersecting 
transversally, and one checks that 
$\omega_Z \cong \cO_Z(-\partial Z)$. 

Since $X_w^v$ has rational singularities, it suffices to show that
$H^j(Z,\pi^*\cL_\lambda)=0$ for $j\ge 1$. By Lemma \ref{ggample}
and the fact that $Z$ is a subvariety of $Z_{\uw}\times Z^{\uv}$,
the boundary $\partial Z$ is the support of an effective ample divisor
$E$ on $Z$. Applying the Kawamata-Viehweg theorem with 
$D:=N\partial Z - E$, where $N$ is a large integer, and
$\cL:=(\pi^*\cL_\lambda)(\partial Z)$, we obtain the desired vanishing
as in the proof of Theorem \ref{bwb}. 

(iii) is checked similarly : let now $E$ be the pull-back on
$Z$ of an effective ample divisor on $Z_{\uw}$ with support 
$\partial Z_{\uw}$. Let $N$ be a large integer, and put
$\cL:=(\pi^*\cL_\lambda)((\partial Z_{\uw})^{\uv})$.
Then the assumptions of the Kawamata-Viehweg theorem are still
verified, since the projection $Z \to Z_{\uw}$ is generically
injective. Thus, we obtain 
$$
H^j(Z,(\pi^*\cL_\lambda)(-(\partial Z^{\uv})_{\uw}))=0 
\quad \text{for} \quad j\ge 1.
$$ 
This implies in turn that 
$$
R^j\pi_*\cO_Z(-(\partial Z^{\uv})_{\uw})=0 
\quad \text{for} \quad j\ge 1.
$$
Together with the isomorphism 
$$
\pi_*\cO_Z(-(\partial Z^{\uv})_{\uw})=\cO_{X_w^v}(-(\partial X^v)_w)
$$
and a Leray spectral sequence argument, this completes the proof.

Likewise, (iv) follows from the vanishing of 
$H^j(Z,(\pi^*L_\lambda) \otimes \omega_Z)$ for $j\ge 1$. In turn, 
this is a consequence of the Grauert-Riemenschneider theorem, since
$L_\lambda$ is ample on $X_w^v$.
\end{proof}

\begin{remarks}
1) One may also show that the restriction 
$H^0(\lambda) \to H^0(X_w^v,L_\lambda)$ is surjective for any dominant
weight $\lambda$. As in Corollary \ref{projnorm}, it follows that the
affine cone over $X_w^v$ has rational singularities in the projective
embedding given by any ample line bundle on $X$. In particular,
$X_w^v$ is projectively normal in any such embedding.

\medskip

\noindent
2) Theorem \ref{ri} (iv) does not extend to all the dominant weights
$\lambda$. Indeed, for $\lambda=0$ we obtain
$$
H^j(X_w^v,\cO_{X_w^v}(-\partial X_w^v)) = H^j(X_w^v,\omega_{X_w^v}).
$$
By Serre duality, this equals 
$H^{\ell(w) - \ell(v) - j}(X_w^v,\cO_{X_w^v})$; i.e., $\bC$ if
$j=\ell(w)-\ell(v)$, and $0$ otherwise, by Theorem \ref{ri} (iii).
\end{remarks}

Next we adapt the construction of Section 3 to obtain a degeneration
of the diagonal of any Richardson variety $X_w^v$. Let
$\lambda:\bC^*\to T$ be as in Subsection 3.1 and let $\cX_w^v$ be the
closure in $X\times X\times \bP^1$ of the subset
$$
\{(x,\lambda(t)x,t) ~\vert~ x\in X_w^v, \;t\in\bC^*\}
\subseteq X\times X\times \bC^*.
$$
We still denote by $\pi:\cX_w^v\to \bP^1$ the projection, then 
$\pi^{-1}(\bC^*)$ identifies again to the product
$\diag(X_w^v) \times \bC^*$ above $\bC^*$. Further, we have the
following analogues of Theorem \ref{fibers} and Proposition \ref{filt}.

\begin{proposition}\label{de}
(i) With the preceding notation, we have equalities of subschemes of
$X\times X$: 
$$
\pi^{-1}(0)=\bigcup_{x\in W,\; v\le x \le w} X_x^v\times X^x_w
\quad \text{and} \quad
\pi^{-1}(\infty)=
\bigcup_{x\in W,\; v\le x \le w} X^x_w\times X_x^v.
$$

\noindent
(ii) The sheaf $\cO_{\pi^{-1}(0)}$ admits a filtration with
associated graded
$$
\bigoplus_{x\in W,\; v\le x \le w} \cO_{X_x^v}\otimes\cO_{X^x_w}
(-(\partial X^x)_w).
$$
Therefore, we have in $K(X\times X)$:
$$
[\cO_{\diag(X_w^v)}] = \sum_{x\in W,\; v\le x \le w}
[\cO_{X_x^v}] \times [\cO_{X^x_w}(-(\partial X^x)_w)].
$$
\end{proposition}

\begin{proof}
Put 
$$
Y_w^v:=\bigcup_{x\in W} X_x^v\times X^x_w.
$$
By the argument of the proof of Theorem \ref{fibers}, we obtain the
inclusion $Y_w^v \subseteq \pi^{-1}(0)$. Further, the proof of
Proposition \ref{filt} shows that the structure sheaf $\cO_{Y_w^v}$
admits a filtration with associated graded given by (ii).

On the other hand, Lemma \ref{const} implies the equality 
$[\cO_{\pi^{-1}(0)}] = [\cO_{\diag(X_w^v)}]$ in $K(X\times X)$.
Further, we have
$$
[\cO_{\diag(X_w^v)}] = [\cO_{\diag(X)}]\cdot [\cO_{X^v\times X_w}]
$$
by Lemma \ref{fp}, since $\diag(X)$ and $X^v\times X_w$ meet properly
in $X\times X$ along $\diag(X_w^v)$. Together with Theorem \ref{dualk}
and Lemma \ref{fp} again, this yields
$$ 
[\cO_{\diag(X_w^v)}] = 
\sum_{x\in W} [\cO_{X_x^v}] \times [\cO_{X_w^x}(-(\partial X^x)_w)]
= [\cO_{Y_w^v}].
$$
Thus, the structure sheaves of $Y_w^v$ and of $\pi^{-1}(0)$ have the
same class in $K(X\times X)$. But we have an exact sequence
$$
0 \to \cF \to \cO_{\pi^{-1}(0)} \to \cO_{Y_w^v} \to 0,
$$
where $\cF$ is a coherent sheaf on $X\times X$. So $[\cF]=0$ in 
$K(X \times X)$, and it follows that $\cF=0$ (e.g., by Example
\ref{db}.1). In other words, $Y_w^v = \pi^{-1}(0)$. This proves (ii)
and the first assertion of (i); the second assertion follows by
symmetry. 
\end{proof}

\subsection{Structure constants and bases of the Grothendieck group}

Let $c_{vw}^x$ be the structure constants of the Grothendieck ring
$K(X)$ in its basis $\{ \cO_w \}$ of Schubert classes, that is, we
have in $K(X)$:
$$
\cO_v \cdot \cO_w=\sum_{x\in W} c_{vw}^x \; \cO_x.
$$
Then Theorem \ref{ri} (i) yields the equality 
$c_{vw}^x = c^x(X_w^{w_o v})$.
Together with Theorem \ref{va}, this implies a solution to Buch's
conjecture: 

\begin{theorem}\label{si}
The structure constants $c_{vw}^x$ satisfy
$$
(-1)^{\ell(v)+\ell(w)+\ell(x)+\ell(w_o)} \; c_{vw}^x \geq 0.
$$
\end{theorem}

Another consequence of Theorem \ref{ri} is the following relation
between the bases $\{\cO_w\}$ and $\{\cI_w\}$ of the group $K(X)$
introduced in 3.4.

\begin{proposition}\label{oi}
We have in $K(X)$
$$
\cO_w = \sum_{v\in W,\; v\le w} \cI_v \quad \text{and} \quad
\cI_w = \sum_{v\in W,\; v\le w} (-1)^{\ell(w)-\ell(v)} \; \cO_v.
$$
\end{proposition}

\begin{proof}
By Theorem \ref{dualk}, we have
$$
\cO_w = \sum_{v\in W} \chi(\cO_w \cdot \cO^v) \; \cI_v.
$$
Further, 
$$
\chi(\cO_w\cdot \cO^v)= \chi(\cO_{X_w^v}) = 
\sum_j (-1)^j \; h^j(\cO_{X_w^v})
$$ 
equals $1$ if $v=w$, and $0$ otherwise, by Theorem \ref{ri}.

Likewise, we obtain 
$$
\cI_w = \sum_{v\in W} \chi(\cI_w \cdot \cI^v) \; \cO_v
\quad \text{and} \quad
\chi(\cI_w \cdot \cI^v) = 
\chi(\cO_{X_w^v}(-\partial X_w^v)) = \chi(\omega_{X_w^v})
$$
by using the equalities 
$\cI_w =[ \cO_{X_w}] - [\cO_{\partial X_w}]$,
$\cI^v =[ \cO_{X^v}] - [\cO_{\partial X^v}]$,
together with Lemma \ref{rs} and Cohen-Macaulayness of Schubert
varieties and their boundaries. Further, we have by Serre duality
and Theorem \ref{ri}:
$$
\chi(\omega_{X_w^v}) = (-1)^{\dim(X_w^v)}\; \chi(\cO_{X_w^v})
= (-1)^{\ell(w) - \ell(v)}.
$$
\end{proof}

\begin{remark}
The preceding proposition implies that the M\"obius function of the
Bruhat order on $W$ maps $(v,w)\in W\times W$ to
$(-1)^{\ell(w)-\ell(v)}$ if $v\le w$, and to $0$ otherwise. We refer to
\cite{Deo77} for a direct proof of this combinatorial fact.
\end{remark}

Using the duality involution $\alpha\mapsto \alpha^{\vee}$ of $K(X)$,
we now introduce another natural basis of this group for which the
structure constants become positive.

\begin{proposition}\label{du}
(i) We have in $K(X)$
$$
[L_\rho\vert_{X_w}(-\partial X_w)] 
= (-1)^{\ell(w_o)-\ell(w)} \; \cO_w^{\vee}.
$$
In particular, the classes 
$$
\cI_w(\rho) := [L_\rho\vert_{X_w}(-\partial X_w)]
= [L_\rho] \cdot \cI_w 
$$ 
form a basis of the Grothendieck group $K(X)$. 

\noindent
(ii)
For any Cohen-Macaulay subscheme $Y$ of $X$ with relative dualizing
sheaf $\omega_{Y/X} = \omega_Y \otimes \omega_X^{-1}$, we have
$$
[\omega_{Y/X}] = \sum_{w\in W} (-1)^{\dim(Y)-\ell(w)} \; 
c^w(Y) \; \cI_w(\rho).
$$
Thus, if $Y$ is a variety with rational singularities, then 
the coordinates of $\omega_{Y/X}$ in the basis $\{\cI_w(\rho)\}$
are the absolute values of the $c^w(Y)$.

\noindent
(iii) The structure constants of $K(X)$ in the basis $\{\cI_w(\rho)\}$
are the absolute values of the structure constants $c_{vw}^x$. 
\end{proposition}

\begin{proof} 
We obtain
$$\displaylines{
[\cO_{X_w}]^{\vee} = (-1)^{\codim(X_w)} \;
[\omega_{X_w}] \cdot [\omega_X^{-1}] 
\hfill\cr\hfill
= (-1)^{\codim(X_w)} \;
[L_{-\rho}\vert_{X_w}(-\partial X_w)]\cdot [L_{2\rho}]
= (-1)^{\ell(w_o)-\ell(w)} \; \cI_w(\rho).
\cr}$$
This proves (i). The assertions (ii), (iii) follow by applying the 
duality involution to Theorems \ref{va} and \ref{si}.
\end{proof}

By similar arguments, we obtain the following relations between the
bases $\{\cI_w(\rho)\}$ and $\{\cO_w\}$.

\begin{proposition}\label{tr}
$$
\cI_w(\rho) = \sum_{v\in W} h_w^v \; \cO_v, 
\quad \text{where} \quad
h_w^v := h^0(X_w^v,L_\rho(-\partial X_w^v)).
$$
In particular, $h_w^v \ne 0$ only if $v\le w$. Further,
$$
\cO_w = \sum_{v\in W,\; v\le w} 
(-1)^{\ell(w)-\ell(v)} h_w^v \; \cI_v(\rho).
$$
\end{proposition}

Next we consider the decomposition of the products 
$[L_\lambda]\cdot \cO_w$ in the basis $\{\cO_v\}$, where $\lambda$ is
a dominant weight. These products also determine the multiplication in
$K(X)$. Indeed, by \cite{Mar76}, this ring is generated by the classes of
line bundles (thus, going to the associated graded 
$\Gr K(X)\cong H^*(X)$, it follows that the cohomology ring is
generated over the rationals by the Chern classes of line bundles).
Since any weight is the difference of two dominant weights, it follows
that the ring $K(X)$ is generated by the classes $[L_\lambda]$, where
$\lambda$ is dominant. This motivates the following:

\begin{theorem}\label{ch}
For any dominant weight $\lambda$ and any $w\in W$, we have in $K(X)$
$$
[L_\lambda] \cdot \cO_w = [L_\lambda\vert_{X_w}] =
\sum_{v\in W,\;v\le w} h^0(X_w^v,L_\lambda(-(\partial X^v)_w)) \; \cO_v.
$$
In particular, the coefficients of $[L_\lambda] \cdot \cO_w$ in the
basis of Schubert classes are non-negative.
\end{theorem}

\begin{proof}
By Theorem \ref{dualk}, we have 
$$
[L_{\lambda}]\cdot\cO_w = \sum_{v\in W} 
\chi([L_{\lambda}] \cdot \cO_w \cdot \cI^v)\; \cO_v.
$$
Further, as in the proof of Proposition \ref{oi}, we obtain
$$
\chi([L_{\lambda}] \cdot \cO_w \cdot \cI^v)
=\chi(X_w^v,L_{\lambda}(-(\partial X^v)_w)).
$$
The latter equals $h^0(X_w^v,L_{\lambda}(-(\partial X^v)_w))$ by
Theorem \ref{ri}.
\end{proof}

Next let $\sigma$ be a non-zero section of $L_\lambda$ on $X$. Then
the structure sheaf of the zero subscheme $Z(\sigma)\subset X$ 
fits into an exact sequence
$$
0 \to L_{-\lambda} \to \cO_X \to \cO_{Z(\sigma)} \to 0.
$$
Thus, the class $[\cO_{Z(\sigma)}] = 1- [L_{-\lambda}]$ depends only
on $\lambda$; we denote this class by $\cO_\lambda$. Note that the
image of $\cO_\lambda$ in the associated graded 
$\Gr K(X) \cong H^*(X)$ is the class of the divisor of $\sigma$, i.e.,
the Chern class $c_1(L_\lambda)$. We now decompose the products 
$\cO_\lambda \cdot \cO_w$ in the basis of Schubert classes.

\begin{proposition}\label{pi}
For any dominant weight $\lambda$ and any $w\in W$, we have in $K(X)$
$$
\cO_\lambda \cdot \cO_w = \sum_{v\in W,\; v<w}
(-1)^{\ell(w)-\ell(v)-1} \; h^0(X_w^v,L_\lambda(-(\partial X_w)^v)) 
\; \cO_v.
$$
\end{proposition}

\begin{proof}
We begin by decomposing the product $[L_\lambda] \cdot \cI_w$ in the
basis $\{\cI_v\}$. As in the proof of Theorem \ref{ch}, we obtain
$$\displaylines{
[L_\lambda] \cdot \cI_w = \sum_{v\in W} 
\chi([L_\lambda] \cdot \cI_w \cdot \cO^v) \; \cI_v 
\hfill\cr\hfill
= \sum_{v\in W,\; v\le w} 
\chi(X_w^v,L_\lambda(-(\partial X_w)^v)) \; \cI_v 
= \sum_{v\in W,\; v\le w} h^0(X_w^v,L_\lambda(-(\partial X_w)^v))
\; \cI_v.
\cr}$$
Applying the duality involution and using the equality 
$$
\cI_w^{\vee} = (-1)^{\ell(w_o)-\ell(w)} \; [L_\rho]\cdot \cO_w,
$$
we obtain
$$
[L_{-\lambda}] \cdot \cO_w = \sum_{v\in W, v\le w}
(-1)^{\ell(w)-\ell(v)} \; h^0(X_w^v,L_\lambda(-(\partial X_w)^v))
\; \cO_v.
$$
Further, $[L_{-\lambda}] = 1 - \cO_\lambda$. Substituting in the 
previous equality completes the proof.
\end{proof}

\begin{remarks}
1) In the case of a fundamental weight $\chi_d$, the divisor of 
the section $p_{w_o \chi_d}$ equals $[X_{w_o s_d}]$, and hence 
$\cO_{w_o \chi_d}$ is the class of the Schubert divisor $X_{w_o s_d}$. 
Thus, Proposition \ref{pi} expresses the structure constants arising 
from the product of the classes of Schubert divisors by arbitrary 
Schubert classes. These structure constants have alternating signs
as predicted by Theorem \ref{si}.

\medskip

\noindent
2) Proposition \ref{pi} gives back the Chevalley formula in 
$H^*(X)$ obtained in Proposition \ref{weights}. 
Indeed, going to the associated graded $\Gr K(X) \cong H^*(X)$ yields
$$
c_1(L_\lambda) \cup [X_w] = 
\sum_v h^0(X_w^v,L_\lambda(-(\partial X_w)^v)) \; [X_v],
$$
the sum over the $v\in W$ such that $v\le w$ and 
$\ell(v) = \ell(w) - 1$. For any such $v$, we know that the Richardson
variety $X_w^v$ is isomorphic to $\bP^1$, identifying the
restriction of $L_\lambda$ to $\cO_{\bP^1}(\lambda_i - \lambda_j)$,
where $v = w s_{ij}$ and $i<j$. Further, $(\partial X_w)^v$
is just the point $vF$, so that 
$L_\lambda\vert_{X_w^v}(-(\partial X_w)^v)$ 
identifies to $\cO_{\bP^1}(\lambda_i - \lambda_j -1)$. Thus, 
$h^0(X_w^v,L_\lambda(-(\partial X_w)^v)) = \lambda_i - \lambda_j$.

\medskip

\noindent
3) The results of this subsection adapt to any partial flag variety
$X = G/P$. In particular, if $X$ is the Grassmannian $\Gr(d,n)$ 
and $L$ is the ample generator of $\Pic(X)$, then we have 
$L_{X_I}(-\partial X_I) \cong \cO_{X_I}$, so that Theorem \ref{ch}
yields the very simple formula
$$
[L\vert_{X_I}] = [L] \cdot \cO_I = \sum_{J,\;J\le I} \cO_J.
$$
In particular, $[L] = \sum_I \cO_I$ (sum over all the multi-indices
$I$). 

By M\"obius inversion, it follows that
$[L^{-1}]\cdot \cO_I = \sum_{J,\;J\le I} \;
(-1)^{\vert I \vert - \vert J\vert} \; \cO_J$. This yields
$$
\cO_{\omega_d} \cdot \cO_I = \sum_{J,\;J < I} 
(-1)^{\vert I \vert - \vert J\vert -1} \; \cO_J,
$$ 
where $\cO_{\omega_d}$ is the class of the Schubert divisor.
\end{remarks}

\medskip

\noindent
{\bf Notes.}
A general reference for this section is \cite{Br02}, from which much of
the exposition is taken.

Stronger versions of Theorem \ref{ri} were obtained in \cite{BrLa03} by
the techniques of Frobenius splitting, and Proposition \ref{de} was
also proved there. These results also follow from standard monomial
theory by work of Kreiman and Lakshmibai for Grassmannians
\cite{KrLa03}, Lakshmibai and Littelmann in general \cite{LaLi03}.

Propositions \ref{oi}, \ref{du} (i) and \ref{tr} are due to Kostant
and Kumar \cite{KoKu90} in the framework of $T$-equivariant $K$-theory;
again, the present approach is also valid in this framework. 

Theorem \ref{ch} also extends readily to $T$-equivariant $K$-theory.
In this form, it is due to Fulton and Lascoux \cite{FuLa94} in the
case of the general linear group. Then the general case was settled by
Pittie and Ram \cite{PiRa99}, Mathieu \cite{Mat00}, Littelmann and
Seshadri \cite{LiSe03}, via very different methods. The latter authors
obtained a more precise version by using standard monomial
theory. Specifically, they constructed a $B$-stable filtration of the 
sheaf $L_\lambda\vert_{X_w}(-\partial X_w)$ with associated graded
sheaf being the direct sum of structure sheaves of Schubert varieties
(with twists by characters). This was generalized to Richardson
varieties by Lakshmibai and Littelmann \cite{LaLi03}, again by using
standard monomial theory. 

This theory constructs bases for spaces of sections of line bundles
over flag varieties, consisting of $T$-eigenvectors which satisfy
very strong compatibility properties to Schubert and opposite Schubert
varieties. It was completed by Littelmann \cite{Li98a}, \cite{Li98b}
after contributions of Lakshmibai, Musili, and Seshadri. Littelmann's
approach is based on methods from combinatorics (the path model in
representation theory) and algebra (quantum groups at roots of
unity). It would be highly desirable to obtain a completely geometric
derivation of standard monomial theory; some steps in this direction
are taken in \cite{BrLa03}.

Another open problem is to obtain a positivity result for the
structure constants of the $T$-equivariant Grothendieck ring. Such a
result would imply both Theorem \ref{si} and Graham's positivity 
theorem \cite{Gra01} for the structure constants in the $T$-equivariant
cohomology ring. A precise conjecture in this direction is formulated
in \cite{GrRa04}, where a combinatorial approach to $T$-equivariant
$K$-theory of flag varieties is developed.

\vfill\eject

\end{document}